\documentclass[12pt,draft,leqno]{article}
\usepackage{amssymb, eucal, latexsym}

0in \oddsidemargin 1,2cm \evensidemargin 0in

\newtheorem{theorem}{Theorem}[section]
\newtheorem{lm}[theorem]{Lemma}
\newtheorem{exa}[theorem]{Example}

\newtheorem{cor}[theorem]{Corollary}
\newtheorem{pro}[theorem]{Proposition}
\newtheorem{defi}[theorem]{Definition}

\newtheorem{nota}[theorem]{Notation}
\newtheorem{notas}[theorem]{Notations}
\newtheorem{rem}[theorem]{Remark}

\newtheorem{fact}[theorem]{Fact}

\newtheorem{nist}[theorem]{}

\def\p{\varphi}
\def\a{\alpha}

\def\d{\delta}

\def\g{\gamma}
\def\GA{\Gamma}

\def\l{\lambda}
\def\LAM{\Lambda}

\def\s{\sigma}

\def\OM{\Omega}

\def\pl{\varphi_\Lambda}

\def\lag{\lambda_A^g}
\def\lbg{\lambda_B^g}

\def\lra{\longrightarrow}

\def\sbe{\subseteq}
\def\spe{\supseteq}
\def\stm{\setminus}
\def\ems{\emptyset}
\def\nes{\neq\emptyset}

\def\cuk{\,\check{}\,}
\def\gek{\,\tilde{}\,}

\def\ex{\exists}
\def\fa{\forall}

\def\we{\wedge}

\def\bv{\bigvee}
\def\cd{\diamond}

\def\ap{^\prime}
\def\inv{^{-1}}
\def\st{\ |\ }

\def\ix{\infty_X}

\def\llx{\ll_{\rho}}
\def\lle{\ll_{\eta}}

\def\llcr{\ll_{C_\rho}}
\def\llce{\ll_{C_\eta}}

\def\nin{\not\in}

\def\card #1{\vert #1 \vert}

\def\ion{i=1,\ldots,n}

\def\CC{{\cal C}}

\def\OO{{\cal O}}

\def\TT{{\cal T}}

\def\Bo{{\bf Bool}}
\def\HC{{\bf HC}}

\def\HLC{{\bf HLC}}

\def\DLC{{\bf DLC}}
\def\CHLC{{\bf CHLC}}
\def\CDLC{{\bf CDLC}}
\def\PDLC{{\bf PDLC}}
\def\DHC{{\bf DHC}}

\def\VAL{{\bf VAL}}
\def\NAL{{\bf NAL}}
\def\PAL{{\bf PAL}}

\def\PLC{{\bf PLC}}

\def\SAC{{\bf DSkeC}}

\def\SKAL{{\bf DSkeLC}}
\def\SHC{{\bf SkeC}}
\def\SKLC{{\bf SkeLC}}

\def\VAL{{\bf DVAL}}

\def\2{\mbox{{\bf 2}}}
\def\3{\mbox{{\bf 3}}}

\def\int{\mbox{{\rm int}}}

\def\cl{\mbox{{\rm cl}}}
\def\CL{\mbox{{\rm Clust}}}
\def\BClu{\mbox{{\rm BClust}}}

\def\ssi{S_{\s\ap}}
\def\vsi{V_{\s\ap}}
\def\jsi{J_{\s\ap}}

\def\doc{\hspace{-1cm}{\em Proof.}~~}
\def\sq{\hspace*{\fill} \hbox{\vrule\vbox{\hrule\phantom{o}\hrule}\vrule}}
\def\sqs{\sq \vspace{2mm}}

\def\BBBB{{\rm I}\!{\rm B}}

\def\eset{\emptyset}   %empty set%
\def\neset{\neq\emptyset} %%not=emptyset%%

\title{{\LARGE\bf
A De Vries-type Duality Theorem for Locally Compact Spaces -- I
}\thanks{This paper was supported by the project no. 136/2008
$``$General and Computer Topology" of the Sofia University $``$St.
Kl. Ohridski".}\\ \vspace{0.35cm}
{\large\bf Georgi D. Dimov}}

\author{}

\date{}

\begin{document}
\maketitle
\begin{abstract}
{\footnotesize
\noindent  A duality theorem for the category of locally compact
Hausdorff spaces and continuous maps which generalizes
 the well-known Duality Theorem of de Vries is proved.}
\end{abstract}

{\footnotesize {\em 2000 MSC:} primary 18A40, 54D45; secondary
06E15, 54C10, 54E05, 06E10.

}

\footnotetext[1]{{\footnotesize {\em E-mail address:}
gdimov@fmi.uni-sofia.bg}}

\baselineskip = \normalbaselineskip

\section*{Introduction}

It is well known that the restriction of the famous localic
duality to the ca\-tegory of locally compact (respectively,
compact) regular frames and frames homomorphisms produces a
duality with the category  $\HLC$ (respectively, $\HC$) of locally
compact (respectively, compact) Hausdorff spaces and continuous
maps (see, e.g., \cite{J}). Another duality for the category
$\HC$, which is in the spirit of the celebrated   Stone Duality
Theorem (\cite{ST}) (namely, that the category of
zero-dimensi\-o\-nal compact Hausdorff spaces and continuous maps
is dually equivalent to the category $\Bo$ of Boolean algebras and
Boolean homomorphisms), was proved by de Vries \cite{dV}. He
introduced the notion of {\em compingent Boolean algebra}\/ and
proved that
 the category $\VAL$ of complete compingent Boolean
algebras and suitable morphisms between them is dually equivalent
to the category $\HC$; thus, the  category $\VAL$ is equivalent,
by the celebrated Gelfand Duality Theorem (\cite{G1,G2,GN,GS}), to
the ca\-tegory of C${}^*$-algebras and algebra homomorphisms. It
is natural to try to extend de Vries Duality Theorem to the
category $\HLC$. An important step in this direction was done by
Roeper \cite{R}. He defined the notion of {\em region-based
topology}\/  as one of the possible formalizations of the ideas of
De Laguna \cite{dL} and Whitehead \cite{W} for a region-based
theory of space.  Following \cite{VDDB, DV}, the region-based
topologies of Roeper appear here as {\em local contact algebras}\/
(briefly, LCAs), because the axioms which they satisfy almost
coincide with the axioms of local proxi\-mities of Leader
\cite{LE}. In his paper \cite{R}, Roeper proved the following
theorem: there is a bijective correspondence between all (up to
homeomorphism) locally compact Hausdorff spaces and all (up to
isomorphism) complete LCAs.  In \cite{D}, using Roeper's Theorem,
a duality theorem for the category $\PLC$ of locally compact
Hausdorff spaces and perfect maps between them was proved. It
implies easily the Duality Theorem of H. de Vries. In the first
part of the present paper, an extension to the category $\HLC$ of
locally compact Hausdorff spaces and continuous maps of both
duality theorems, obtained in \cite{dV} and \cite{D}, is
presented. A duality theorem for the category of connected locally
compact Hausdorff spaces and continuous maps is proved as well. In
the second part \cite{D5} of the paper, some applications of the
methods and results of the first part as well as of \cite{ST, dV,
R, D,D1,D2} are given. In particular, a slight generalization of
the Stone's results from \cite{ST} concerning the extension of his
Duality Theorem to the category of locally compact
zero-dimensional Hausdorff spaces is obtained; a description of
the products in the category $\DLC$ which is dual to the category
$\HLC$ is presented; a completion theorem for
 LCAs is proved; a characterization of the dual objects of
 the metrizable locally compact spaces is found;
  a direct proof of  Ponomarev's solution \cite{P1} of
  Problem 72 of  G. Birkhoff (see \cite{Bi}) is obtained, and the class of spaces which
 are
 co-absolute with the (zero-dimensional) Eberlein compacts is
 described.

We now fix the notations. They will be used in both parts of the
paper.

If $\CC$ denotes a category, we write $X\in \card\CC$ if $X$ is
 an object of $\CC$, and $f\in \CC(X,Y)$ if $f$ is a morphism of
 $\CC$ with domain $X$ and codomain $Y$.

All lattices are with top (= unit) and bottom (= zero) elements,
denoted respectively by 1 and 0. We do not require the elements
$0$ and $1$ to be distinct. We set $\2=\{0,1\}$, where $0\neq 1$.
If $(A,\le)$ is a poset and $a\in A$, we set
$\downarrow_A(a)=\{b\in A\st b\le a\}$ (we will write even
$``\downarrow(a)$" instead of $``\downarrow_A(a)$" when there is
no ambiguity); if $B\sbe A$ then we set
$\downarrow(B)=\bigcup\{\downarrow(b)\st b\in B\}$.

If $X$ is a set then we denote the power set of $X$ by $P(X)$. If
$Y$ is also a set and $f:X\lra Y$ is a function, then we will set,
for every $U\sbe X$,
 $f^\sharp(U)=\{y\in Y\st f\inv(y)\sbe U\}$.
 If
$(X,\tau)$ is a topological space and $M$ is a subset of $X$, we
denote by $\cl_{(X,\tau)}(M)$ (or simply by $\cl(M)$ or
$\cl_X(M)$) the closure of $M$ in $(X,\tau)$ and by
$\int_{(X,\tau)}(M)$ (or briefly by $\int(M)$ or $\int_X(M)$) the
interior of $M$ in $(X,\tau)$. The Alexandroff compactification of
a locally compact Hausdorff non-compact space $X$ will be denoted
by $\a X$ and the added point by $\ix$ (i.e. $\a X=X\cup\{\ix\}$).
The (positive) natural numbers are denoted by $\mathbb{N}$
($\mathbb{N}^+$) and the real line -- by $\mathbb{R}$.

The  closed maps  between topological spaces are assumed to be
continuous but are not assumed to be onto. Recall that a map is
{\em perfect}\/ if it is  compact (i.e. point inverses are compact
sets) and closed. A continuous map $f:X\lra Y$ is {\em
irreducible}\/ if $f(X)=Y$ and for each proper closed subset $A$
of $X$, $f(A)\neq Y$.

\section{Preliminaries}
%%%%%%%%%%%%%%%%%%%%%%%%%%%%%%%%%%%%%%%%%%%%%%%%%%%%%%%%%%%%
%%%% Section 1. %%%%%%%%%%%%%%%%%%%%%%%
%%%%%%%%%%%%%%%%%%%%%%%%%%%%%%%%%%%%%%%%%%%%%%%%%%%%%%%%%%%%
%

\begin{defi}\label{conalg}
\rm
An algebraic system $\underline {B}=(B,0,1,\vee,\we, {}^*, C)$ is
called a {\it contact Boolean algebra}\/ or, briefly, {\it contact
algebra} (abbreviated as CA) (\cite{DV})
 if
$(B,0,1,\vee,\we, {}^*)$ is a Boolean algebra (where the operation
$``$complement" is denoted by $``\ {}^*\ $")
  and $C$
is a binary relation on $B$, satisfying the following axioms:

\smallskip

\noindent (C1) If $a\not= 0$ then $aCa$;\\
(C2) If $aCb$ then $a\not=0$ and $b\not=0$;\\
(C3) $aCb$ implies $bCa$;\\
(C4) $aC(b\vee c)$ iff $aCb$ or $aCc$.

\smallskip

\noindent We shall simply write $(B,C)$ for a contact algebra. The
relation $C$  is called a {\em  contact relation}. When $B$ is a
complete Boolean algebra, we will say that $(B,C)$ is a {\em
complete contact Boolean algebra}\/ or, briefly, {\em complete
contact algebra} (abbreviated as CCA). If $a\in B$ and $D\sbe B$,
we will write $``aCD$" for $``(\fa d\in D)(aCd)$".

We will say that two CA's $(B_1,C_1)$ and $(B_2,C_2)$ are  {\em
CA-isomorphic} iff there exists a Boolean isomorphism $\p:B_1\lra
B_2$ such that, for each $a,b\in B_1$, $aC_1 b$ iff $\p(a)C_2
\p(b)$. Note that in this paper, by a $``$Boolean isomorphism" we
understand an isomorphism in the category $\Bo$.

A CA $(B,C)$  is called {\em connected}\/ if it satisfies the
following axiom:

\smallskip

\noindent (CON) If $a\neq 0,1$ then $aCa^*$.

\smallskip

A contact algebra $(B,C)$ is called a {\it  normal contact Boolean
algebra}\/ or, briefly, {\it  normal contact algebra} (abbreviated
as NCA) (\cite{dV,F}) if it satisfies the following axioms (we
will write $``-C$" for $``not\ C$"):

\smallskip

\noindent (C5) If $a(-C)b$ then $a(-C)c$ and $b(-C)c^*$ for some $c\in B$;\\
(C6) If $a\not= 1$ then there exists $b\not= 0$ such that
$b(-C)a$.

\smallskip

\noindent A normal CA is called a {\em complete normal contact
Boolean algebra}\/ or, briefly, {\em complete normal contact
algebra} (abbreviated as CNCA) if it is a CCA. The notion of
normal contact algebra was introduced by Fedorchuk \cite{F} under
the name {\em Boolean $\d$-algebra}\/ as an equivalent expression
of the notion of compingent Boolean algebra of de Vries. We call
such algebras $``$normal contact algebras" because they form a
subclass of the class of contact algebras and naturally arise in
normal Hausdorff spaces.

Note that if $0\neq 1$ then the axiom (C2) follows from the axioms
(C6) and (C4).

For any CA $(B,C)$, we define a binary relation  $``\ll_C $"  on
$B$ (called {\em non-tangential inclusion})  by $``\ a \ll_C b
\leftrightarrow a(-C)b^*\ $". Sometimes we will write simply
$``\ll$" instead of $``\ll_C$".
\end{defi}

The relations $C$ and $\ll$ are inter-definable. For example,
normal contact algebras could be equivalently defined (and exactly
in this way they were introduced (under the name of {\em
compingent Boolean algebras}) by de Vries in \cite{dV}) as a pair
of a Boolean algebra $B=(B,0,1,\vee,\we,{}^*)$ and a binary
relation $\ll$ on $B$ subject to the following axioms:

\smallskip

\noindent ($\ll$1) $a\ll b$ implies $a\leq b$;\\
($\ll$2) $0\ll 0$;\\
($\ll$3) $a\leq b\ll c\leq t$ implies $a\ll t$;\\
($\ll$4) $a\ll c$ and $b\ll c$ implies $a\vee b\ll c$;\\
($\ll$5) If  $a\ll c$ then $a\ll b\ll c$  for some $b\in B$;\\
($\ll$6) If $a\neq 0$ then there exists $b\neq 0$ such that $b\ll
a$;\\
($\ll$7) $a\ll b$ implies $b^*\ll a^*$.

\smallskip

Note that if $0\neq 1$ then the axiom ($\ll$2) follows from the
axioms ($\ll$3), ($\ll$4), ($\ll$6) and ($\ll$7).

\smallskip

Obviously, contact algebras could be equivalently defined as a
pair of a Boolean algebra $B$ and a binary relation $\ll$ on $B$
subject to the  axioms ($\ll$1)-($\ll$4) and ($\ll$7).

\smallskip

It is easy to see that axiom (C5) (resp., (C6)) can be stated
equivalently in the form of ($\ll$5) (resp., ($\ll$6)).

\begin{exa}\label{extrcr}
\rm Let $B$ be a Boolean algebra. Then there exist a largest and a
smallest contact relations on $B$; the largest one, $\rho_l$, is
defined by $a\rho_l b$ iff $a\neq 0$ and $b\neq 0$, and the
smallest one, $\rho_s$, by $a\rho_s b$ iff $a\wedge b\neq 0$.

Note that, for $a,b\in B$, $a\ll_{\rho_s} b$ iff $a\le b$; hence
$a\ll_{\rho_s} a$, for any $a\in B$. Thus $(B,\rho_s)$ is a normal
contact algebra.
\end{exa}

\begin{exa}\label{rct}
\rm Recall that a subset $F$ of a topological space $(X,\tau)$ is
called {\em regular closed}\/ if $F=\cl(\int (F))$. Clearly, $F$
is regular closed iff it is the closure of an open set.

For any topological space $(X,\tau)$, the collection $RC(X,\tau)$
(we will often write simply $RC(X)$) of all regular closed subsets
of $(X,\tau)$ becomes a complete Boolean algebra
$(RC(X,\tau),0,1,\we,\vee,{}^*)$ under the following operations:
$$ 1 = X,  0 = \emptyset, F^* = \cl(X\stm F), F\vee G=F\cup G,
F\we G =\cl(\int(F\cap G)).
$$
The infinite operations are given by the following formulas:
$\bigvee\{F_\g\st \g\in\GA\}=\cl(\bigcup\{F_\g\st
\g\in\GA\})(=\cl(\bigcup\{\int(F_\g)\st \g\in\GA\})),$ and
$\bigwedge\{F_\g\st \g\in\GA\}=\cl(\int(\bigcap\{F_\g\st
\g\in\GA\})).$

It is easy to see that setting $F \rho_{(X,\tau)} G$ iff $F\cap
G\not = \ems$, we define a contact relation $\rho_{(X,\tau)}$ on
$RC(X,\tau)$; it is called a {\em standard contact relation}. So,
$(RC(X,\tau),\rho_{(X,\tau)})$ is a CCA (it is called a {\em
standard contact algebra}). We will often write simply $\rho_X$
instead of $\rho_{(X,\tau)}$. Note that, for $F,G\in RC(X)$,
$F\ll_{\rho_X}G$ iff $F\sbe\int_X(G)$.

Clearly, if $(X,\tau)$ is a normal Hausdorff space then the
standard contact algebra $(RC(X,\tau),\rho_{(X,\tau)})$ is a
complete NCA.

A subset $U$ of $(X,\tau)$ such that $U=\int(\cl(U))$ is said to
be {\em regular open}. The set of all regular open subsets of
$(X,\tau)$ will be denoted by $RO(X,\tau)$ (or briefly, by
$RO(X)$). Define Boolean operations and contact $\d_X$ in $RO(X)$
as follows: $U\vee V=\int(\cl(U\cup V))$, $U\wedge V=U\cap V$,
$U^{*}=\int(X\setminus U)$, $0=\eset$, $1=X$ and $U\d_X V$ iff
$\cl (U)\cap \cl (V)\neset$. Then $(RO(X),\d_X)$ is a CA. This
algebra is also complete, considering the infinite meet
$\bigwedge\{U_i\st i\in I\} =\int(\bigcap_{i\in I}U_{i})$.

Note that $(RO(X),\d_X)$ and $(RC(X),\rho_X)$ are isomorphic CAs.
The isomorphism $f$ between them is defined by  $f(U)=\cl(U)$, for
every $U\in RO(X)$.
\end{exa}

We will need a lemma from \cite{CNG}:

\begin{lm}\label{isombool}
Let $X$ be a dense subspace of a topological space $Y$. Then the
functions $r:RC(Y)\lra RC(X)$, $F\mapsto F\cap X$, and
$e:RC(X)\lra RC(Y)$, $G\mapsto \cl_Y(G)$, are Boolean isomorphisms
between Boolean algebras $RC(X)$ and $RC(Y)$, and $e\circ
r=id_{RC(Y)}$, $r\circ e=id_{RC(X)}$.
\end{lm}

\begin{fact}\label{confact}{\rm (\cite{BG})}
Let $(X,\tau)$ be a topological space. Then the standard contact
algebra $(RC(X,\tau),\rho_{(X,\tau)})$ is connected iff the space
$(X,\tau)$ is connected.
\end{fact}

The following notion is a lattice-theoretical counterpart of the
corresponding notion from the theory of proximity spaces (see
\cite{NW}):

\begin{nist}\label{defcluclan}
\rm Let $(B,C)$ be a CA. Then  a non-empty subset $\s $ of $B$ is
called a {\em cluster in} $(B,C)$
%(see \cite{VDDB})
if the
following conditions are satisfied:

\smallskip

\noindent (K1) If $a,b\in\s $ then $aCb$;\\
(K2) If $a\vee b\in\s $ then $a\in\s $ or $b\in\s $;\\
(K3) If $aCb$ for every $b\in\s $, then $a\in\s $.

\smallskip

\noindent The set of all clusters in $(B,C)$ will be denoted  by
$\CL(B,C)$.
\end{nist}

The next  assertion can be proved exactly as Lemma 5.6 of
\cite{NW}:

\begin{fact}\label{fact29}
%{\rm (\cite{VDDB})}
If $\s_1,\s_2$ are two clusters in a normal contact algebra
$(B,C)$ and $\s_1\sbe \s_2$ then $\s_1=\s_2$.
\end{fact}

\begin{theorem}\label{conclustth}{\rm (\cite{VDDB})}
A subset $\s$ of a normal contact algebra $(B,C)$ is a cluster iff
there exists an ultrafilter $u$ in $B$ such that
\begin{equation}\label{ultclu}
\s=\{a\in B\st aCb \mbox{ for every } b\in u\}.
\end{equation}
Moreover, given $\s$ and $a_0\in \s$, there exists an ultrafilter
$u$ in $B$ satisfying (\ref{ultclu}) which contains $a_0$.

Note that everywhere in this assertion,  the word $``$ultrafilter"
 can be replaced by $``$basis of an ultrafilter".
\end{theorem}

\begin{cor}\label{uniqult}{\rm (\cite{VDDB})}
Let $(B,C)$ be a normal contact algebra and $u$ be an ultrafilter
(or a basis of an ultrafilter) in $B$. Then there exists a unique
cluster $\s_u$ in $(B,C)$ containing $u$, and
\begin{equation}\label{sigmau}
\s_u=\{a\in B\st aCb \mbox{  for every } b\in u\}.
\end{equation}
\end{cor}

\begin{pro}\label{cluendcor}{\rm (\cite{R,D1})}
Let $(B,C)$ be a normal contact algebra, $\s$ be a cluster in
$(B,C)$, $a\in B$ and $a\not\in\s$. Then there exists  $b\in B$
such that $b\not\in\s$ and $a\ll b$.
 \end{pro}

The following notion is a lattice-theoretical counterpart of the
Leader's notion of {\em local proximity} (\cite{LE}):

\begin{defi}\label{locono}{\rm (\cite{R})}
\rm An algebraic system $\underline {B}_{\, l}=(B,0,1,\vee,\we,
{}^*, \rho, \BBBB)$ is called a {\it local contact Boolean
algebra}\/ or, briefly, {\it local contact algebra} (abbreviated
as LCA)   if $(B,0,1, \vee,\we, {}^*)$ is a Boolean algebra,
$\rho$ is a binary relation on $B$ such that $(B,\rho)$ is a CA,
and $\BBBB$ is an ideal (possibly non proper) of $B$, satisfying
the following axioms:

\smallskip

\noindent(BC1) If $a\in\BBBB$, $c\in B$ and $a\ll_\rho c$ then
$a\ll_\rho b\ll_\rho c$ for some $b\in\BBBB$  (see \ref{conalg}
for
$``\ll_\rho$");\\
(BC2) If $a\rho b$ then there exists an element $c$ of $\BBBB$
such that
$a\rho (c\we b)$;\\
(BC3) If $a\neq 0$ then there exists  $b\in\BBBB\stm\{0\}$ such
that $b\ll_\rho a$.

\smallskip

We shall simply write  $(B, \rho,\BBBB)$ for a local contact
algebra.  We will say that the elements of $\BBBB$ are {\em
bounded} and the elements of $B\stm \BBBB$  are  {\em unbounded}.
When $B$ is a complete Boolean algebra,  the LCA $(B,\rho,\BBBB)$
is called a {\em complete local contact Boolean algebra}\/ or,
briefly, {\em complete local contact algebra} (abbreviated as
CLCA).

We will say that two local contact algebras $(B,\rho,\BBBB)$ and
$(B_1,\rho_1,\BBBB_1)$ are  {\em LCA-isomorphic} if there exists a
Boolean isomorphism $\p:B\lra B_1$ such that, for $a,b\in B$,
$a\rho b$ iff $\p(a)\rho_1 \p(b)$, and $\p(a)\in\BBBB_1$ iff
$a\in\BBBB$. A map $\p:(B,\rho,\BBBB)\lra(B_1,\rho_1,\BBBB_1)$ is
called an {\em LCA-embedding} if $\p:B\lra B_1$ is an injective
Boolean homomorphism (i.e., Boolean monomorphism) and, moreover,
for any $a,b\in B$, $a\rho b$ iff $\p(a)\rho_1 \p(b)$, and
$\p(a)\in\BBBB_1$ iff $a\in\BBBB$.

An LCA $(B,\rho,\BBBB)$ is called {\em connected}\/ if the CA
$(B,\rho)$ is connected.
\end{defi}

\begin{rem}\label{conaln}
\rm Note that if $(B,\rho,\BBBB)$ is a local contact algebra and
$1\in\BBBB$ then $(B,\rho)$ is a normal contact algebra.
Conversely, any normal contact algebra $(B,C)$ can be regarded as
a local contact algebra of the form $(B,C,B)$.
\end{rem}

The following lemmas
%from \cite{VDDB}
are lattice-theoretical
counterparts of some theorems from Leader's paper \cite{LE}.

\begin{defi}\label{Alexprn}{\rm (\cite{VDDB})}
\rm Let $(B,\rho,\BBBB)$ be a local contact algebra. Define a
binary relation $``C_\rho$" on $B$ by
\begin{equation}\label{crho}
aC_\rho b\ \mbox{ iff }\ a\rho b\ \mbox{ or }\ a,b\not\in\BBBB.
\end{equation}
It is called the\/ {\em Alexandroff extension of}\/ $\rho$ {\em
relatively to the LCA} $(B,\rho,\BBBB)$ (or, when there is no
ambiguity, simply, the \/ {\em Alexandroff extension of}\/
$\rho$).
\end{defi}

\begin{lm}\label{Alexprn1}{\rm (\cite{VDDB})}
Let $(B,\rho,\BBBB)$ be a local contact algebra. Then
$(B,C_\rho)$, where $C_\rho$ is the Alexandroff extension of
$\rho$, is a normal contact algebra.
\end{lm}

\begin{defi}\label{boundcl}
\rm Let $(B,\rho,\BBBB)$ be a local contact algebra. We will say
that $\s$ is a {\em cluster in} $(B,\rho,\BBBB)$ if $\s$ is a
cluster in the NCA $(B,C_\rho)$ (see \ref{Alexprn} and
\ref{Alexprn1}). A cluster $\s$ in $(B,\rho,\BBBB)$ (resp., an
ultrafilter in $B$) is called {\em bounded}\/ if $\s\cap\BBBB\nes$
(resp., $u\cap\BBBB\nes$). The set of all bounded clusters in
$(B,\rho,\BBBB)$ will be denoted by $\BClu(B,\rho,\BBBB)$.
\end{defi}

\begin{lm}\label{neogrn}{\rm (\cite{VDDB})}
Let $(B,\rho,\BBBB)$ be a local contact algebra and let
$1\not\in\BBBB$. Then $\s_\infty^{(B,\rho,\BBBB)}=\{b\in B\st
b\not\in\BBBB\}$ is a cluster in $(B,\rho,\BBBB)$. (Sometimes we
will simply write $\s_\infty$
%or $\s_\infty^B$
instead of $\ \s_\infty^{(B,\rho,\BBBB)}$.)
\end{lm}

\begin{fact}\label{bstar}
Let $(B,\rho,\BBBB)$ be a local contact algebra and  $\s$ be a
bounded cluster in it (see \ref{boundcl}). Then there exists $b\in
\BBBB$  such that  $b^*\nin\s$.
\end{fact}

\doc
Let $b_0\in\s\cap\BBBB$. Since $b_0\ll_\rho 1$, (BC1) implies that
there exists $b\in\BBBB$ such that $b_0\ll_\rho b$. Then
$b_0(-\rho)b^*$ and since $b_0\in\BBBB$, we obtain that
$b_0(-C_\rho)b^*$. Thus $b^*\nin\s$. \sqs

\begin{nota}\label{compregn}
\rm Let $(X,\tau)$ be a topological space. We denote by
$CR(X,\tau)$ the family of all compact regular closed subsets of
$(X,\tau)$. We will often write  $CR(X)$ instead of $CR(X,\tau)$.

 If $x\in X$ then we
set:
\begin{equation}\label{sxvx}
\s_x=\{F\in RC(X)\st x\in F\}.
\end{equation}
\end{nota}

\begin{fact}\label{stanlocn}(\cite{R,VDDB})
Let $(X,\tau)$ be a locally compact Hausdorff space. Then:

(a) the triple
$$(RC(X,\tau),\rho_{(X,\tau)}, CR(X,\tau))$$
 (see \ref{rct} for $\rho_{(X,\tau)}$)
  is a complete local contact algebra; it is called a
{\em standard local contact algebra};

(b) for every $x\in X$, $\s_x$ is a bounded cluster in the
standard local contact algebra $(RC(X,\tau),\rho_{(X,\tau)},
CR(X,\tau))$.
\end{fact}

We will often use (even without citing it explicitly) the
following well-known assertion (see, e.g., \cite[Theorem
3.3.2]{E}):

\begin{pro}\label{loccome}
For every compact subspace $K$ of a locally compact space $X$ and
every open set $V\sbe X$ that contains $K$ there exists an open
set $U\sbe X$ such that $K\sbe U\sbe \cl(U)\sbe V$ and $\cl(U)$ is
compact.
\end{pro}

For all notions and notations not defined here see \cite{AHS, J,
E, NW, Si}.

\section{The Generalization of De Vries Duality Theorem}

The next  theorem was proved by Roeper \cite{R} (but its
particular cases concerning compact Hausdorff spaces and NCAs were
proved by de Vries \cite{dV}).  We will give a sketch of its
proof; it follows the plan of the proof presented in \cite{VDDB}.
The notations and the facts stated here will be used later on.

Recall that if $(A,\le)$ is a poset and $B\sbe A$ then $B$ is said
to be a {\em dense subset of} $A$ if for any $a\in A\stm\{0\}$
there exists $b\in B\stm\{0\}$ such that $b\le a$; when
$(B,\le_1)$ is a poset and $f:A\lra B$ is a map, then we will say
that $f$ is a {\em dense map}\/ if $f(A)$ is a dense subset of
$B$.

\begin{theorem}\label{roeperl}{\rm (P. Roeper \cite{R}
for locally compact spaces and de Vries \cite{dV} for compact
spaces)}

\noindent(a) \ For any LCA (resp., NCA) $(B,\rho,\BBBB)$ there
exists a locally compact (resp., compact) Hausdorff space
$$X=\Psi^a(B,\rho,\BBBB)$$ and a dense LCA-embedding
$\l_{(B,\rho,\BBBB)}:(B,\rho,\BBBB)\lra(RC(X),\rho_X,CR(X))$
(hence, if $(B,\rho,\BBBB)$ is a CLCA then $\l_{(B,\rho,\BBBB)}$
is an LCA-isomorphism).

\smallskip

\noindent(b) \ There exists a bijective correspondence between the
class of all (up to isomorphism) CLCAs  and the class of all (up
to homeomorphism) locally compact Hausdorff spaces; its
restriction to the class of all (up to isomorphism) CNCAs gives a
bijective correspondence between the later class and the class of
all (up to homeomorphism) compact Hausdorff spaces.
\end{theorem}

\noindent{\em Sketch of the Proof.}~ (A) Let $(X,\tau)$ be a
locally compact Hausdorff space. We put
\begin{equation}\label{psit1}
\Psi^t(X,\tau)=(RC(X,\tau),\rho_{(X,\tau)},CR(X,\tau))
\end{equation}
(see \ref{stanlocn} and \ref{compregn} for the notations).

\noindent(B)~ Let $(B,\rho,\BBBB)$ be a
%complete
local contact algebra. Let $C=C_\rho$ be the Alexandroff extension
of $\rho$ (see \ref{Alexprn1}). Then, by  \ref{Alexprn1}, $(B,C)$
is a
%complete
normal contact algebra. Put $X=\CL(B,C)$ and let $\TT$
be the topology on $X$ having as a closed base the family
$\{\l_{(B,C)}(a)\st a\in B\}$ where, for every $a\in B$,
\begin{equation}\label{h}
\l_{(B,C)}(a) = \{\s \in X\st  a \in \s\}.
\end{equation}
Sometimes we will write simply $\l_B$ instead of $\l_{(B,C)}$.

\noindent Note that
\begin{equation}\label{intha}
X\stm \l_B(a)= \int(\l_B(a^*)),
\end{equation}
\begin{equation}\label{ee}
\mbox{the family } \{\int(\l_B(a))\st a\in B\} \mbox{ is an open
base of }(X,\TT)
\end{equation}
and, for every $a\in B$,
\begin{equation}\label{haregcl}
\l_B(a)\in RC(X,\TT).
\end{equation}
It can be proved that
\begin{equation}\label{embedd}
\l_B:(B,C)\lra (RC(X),\rho_X) \mbox{ is a dense CA-embedding}
\end{equation}
and hence,
\begin{equation}\label{isom}
\mbox{when } B \mbox{ is a complete Boolean algebra, } \l_B \mbox{
is a CA-isomorphism.}
\end{equation}
Further,
\begin{equation}\label{xcomp}
(X,\TT) \mbox{ is a compact Hausdorff space.}
\end{equation}

\noindent(B1)~ Let $1\in\BBBB$. Then $C=\rho$ and $\BBBB=B$, so
that $(B,\rho,\BBBB)=(B,C,B)=(B,C)$ is a
%complete
normal contact
algebra (see \ref{conaln}), and we put
\begin{equation}\label{phiapcn}
\Psi^a(B,\rho,\BBBB)(=\Psi^a(B,C,B)=\Psi^a(B,C))=(X,\TT).
\end{equation}

\medskip

\noindent(B2)~ Let $1\not\in\BBBB$. Then, by Lemma \ref{neogrn},
the set $\s_\infty=\{b\in B\st b\not\in\BBBB\}$ is a cluster in
$(B,C)$ and, hence, $\s_\infty\in X$.  Let $L=X\stm\{\s_\infty\}$.
Then
%$L=\BClu(B,\rho,\BBBB)$, i.e.
\begin{equation}\label{L}
  L \mbox{ is the set of all
bounded clusters of } (B,\rho,\BBBB)
\end{equation}
(sometimes we will write $L_{(B,\rho,\BBBB)}$ or $L_B$ instead of
$L$);
 let the topology $\tau(=\tau_{(B,\rho,\BBBB)})$ on $L$ be the
subspace topology, i.e. $\tau=\TT_{|_L} $. Then $(L,\tau)$ is a
locally compact Hausdorff space. We put
\begin{equation}\label {phiapc}
\Psi^a(B,\rho,\BBBB)=(L,\tau).
\end{equation}
Let
\begin{equation}\label{hapni}
\l^l_{(B,\rho,\BBBB)}(a)=\l_{(B,C_\rho)}(a)\cap L,
\end{equation}
for each $a\in B$. We will write simply $\l^l_B$ (or even
$\l_{(B,\rho,\BBBB)}$) instead of $\l^l_{(B,\rho,\BBBB)}$ when
this does not lead to ambiguity. One can show that:

\smallskip

\noindent (I) $L$ is a dense subset of the topological space
$X$;\\
(II) $\l^l_B$ is a dense Boolean monomorphism of the Boolean
algebra $B$ in
the Boolean algebra $RC(L,\tau)$ (and, hence, when $B$ is complete, then $\l^l_B$
is a Boolean isomorphism);\\
(III) $b\in\BBBB$ iff $\l^l_B(b)\in CR(L)$;\\
(IV) $a\rho b$ iff $\l^l_B(a)\cap \l^l_B(b)\neq\ems$.\\
Hence, $X$ is the Alexandroff (i.e. one-point) compactification of
$L$ and
\begin{equation}\label{hapisomem}
\ \ \ \l^l_B: (B,\rho,\BBBB)\lra (RC(L),\rho_L, CR(L)) \mbox{ is a
dense LCA-embedding;}
\end{equation}
thus,
\begin{equation}\label{hapisom}
\mbox{ when } (B,\rho,\BBBB)\mbox{ is a CLCA then } \l^l_B \mbox{
is an LCA-isomorphism.}
\end{equation}
 Note also that for every $b\in B$,
\begin{equation}\label{intl}
\int_{L_B}(\l^l_B(b))=L_B\cap\int_X(\l_B(b)).
\end{equation}
Using (\ref{intha}) and (\ref{intl}), we get readily that for
every $b\in B$,
\begin{equation}\label{inthal}
L\stm \l_B^l(b)= \int_L(\l_B^l(b^*)).
\end{equation}

\medskip

\noindent(C)~ For every LCA $(B,\rho,\BBBB)$ and every $a\in B$,
set
\begin{equation}\label{lbg}
\l^g_{(B,\rho,\BBBB)}(a)=\l_{(B,C_\rho)}(a)\cap\Psi^a(B,\rho,\BBBB).
\end{equation}
We will write simply $\l^g_B$ instead of $\l^g_{(B,\rho,\BBBB)}$
when this does not lead to ambiguity.
 Thus, when $1\in\BBBB$,
we have that $\l^g_B=\l_B$, and  if $1\nin\BBBB$  then
$\l^g_B=\l^l_B$. Hence,
\begin{equation}\label{hapisomnemb}
\ \ \l^g_B: (B,\rho,\BBBB)\lra (\Psi^t\circ\Psi^a)(B,\rho,\BBBB)
\mbox{ is a dense LCA-embedding,}
\end{equation}
and
\begin{equation}\label{hapisomn}
\mbox{when } (B,\rho,\BBBB)\mbox{ is a CLCA then }\l^g_B \mbox{ is
an LCA-isomorphism.}
\end{equation}

With the next assertion we specify (\ref{ee}):
\begin{equation}\label{eel}
 \{\int_{\Psi^a(B,\rho,\BBBB)}(\l_B^g(a))\st a\in \BBBB\} \mbox{
is an open base of } \Psi^a(B,\rho,\BBBB).
\end{equation}
Note that (\ref{embedd}) and (IV)
%(see (B2) above)
imply that if
$(B,\rho,\BBBB)$ is an LCA then
\begin{equation}\label{pbcl}
 a\rho b \mbox{ iff there exists  } \s\in\Psi^a(B,\rho,\BBBB)\mbox{ such
 that } a,b\in\s.
\end{equation}

\noindent(D)~~ Let $(X,\tau)$ be a compact Hausdorff space. Then
it can be proved that the map
\begin{equation}\label{nison}
t_{(X,\tau)}:(X,\tau)\lra\Psi^a(\Psi^t(X,\tau)),
\end{equation}
defined by  $t_{(X,\tau)}(x)=\{F\in RC(X,\tau)\st x\in
F\}(=\s_x)$, for all $x\in X$, is a homeomorphism (we will also
write simply $t_X$ instead of $t_{(X,\tau)}$).

Let  $(L,\tau)$ be a non-compact locally compact Hausdorff space.
Put $B=RC(L,\tau)$, $\BBBB=CR(L,\tau)$ and $\rho=\rho_L$. Then
$(B,\rho,\BBBB)=\Psi^t(L,\tau)$ and $1\nin\BBBB$ (here $1=L$). It
can be shown that the map
\begin{equation}\label{homeo}
t_{(L,\tau)}:(L,\tau)\lra\Psi^a(\Psi^t(L,\tau)),
\end{equation}
defined by  $t_{(L,\tau)}(x)=\{F\in RC(L,\tau)\st x\in
F\}(=\s_x)$, for all $x\in L$, is a homeomorphism; we will often
write simply $t_L$ instead of $t_{(L,\tau)}$.

Therefore $\Psi^a(\Psi^t(L,\tau))$ is homeomorphic to $(L,\tau)$
and when $(B,\rho,\BBBB)$ is a CLCA then
$\Psi^t(\Psi^a(B,\rho,\BBBB))$ is LCA-isomorphic to
$(B,\rho,\BBBB)$.
 \sqs

Note that (\ref{hapisomnemb}) and (\ref{eel}) imply that if
$(B,\rho,\BBBB)$ is an LCA, $X=\Psi^a(B,\rho,\BBBB)$ and
$(A,\eta,\BBBB\ap)=\lbg(B,\rho,\BBBB)$ then for every $a\in
RC(X)$,
\begin{equation}\label{aunion}
a=\bigvee\{b\in\BBBB\ap\st b\ll_{\rho_X} a\}.
\end{equation}
In particular, for every $a\in B$,
\begin{equation}\label{aunionp}
a=\bigvee\{b\in\BBBB\st b\ll_{\rho} a\}.
\end{equation}

In \cite{D2}, a description of the frame of all open subsets of a
locally compact Hausdorff space $X$ in terms of the corresponding
to it CLCA $\Psi^t(X)$ is given (see (\ref{psit1}) for the
notation $\Psi^t$). Since we need this description, we will recall
it. In order to make the paper more self-contained, we will even
give the necessary proofs.

\begin{defi}\label{lideal}{\rm (\cite{D2})}
\rm Let $(A,\rho,\BBBB)$ be an LCA. An ideal $I$ of $A$ is called
a {\em $\d$-ideal} if $I\sbe \BBBB$ and for any $a\in I$ there
exists $b\in I$ such that $a\ll_\rho b$. If $I_1$ and $I_2$ are
two $\d$-ideals of $(A,\rho,\BBBB)$ then we put $I_1\le I_2$ iff
$I_1\sbe I_2$. We will denote by $(I(A,\rho,\BBBB),\le)$ the poset
of all $\d$-ideals of $(A,\rho,\BBBB)$.
\end{defi}

\begin{fact}\label{dideal}{\rm (\cite{D2})}
Let $(A,\rho,\BBBB)$ be an LCA. Then, for every $a\in A$, the set
$I_a=\{b\in\BBBB\st b\ll_\rho a\}$ is a $\d$-ideal. Such
$\d$-ideals will be called\/ {\em principal $\d$-ideals}.
\end{fact}

Recall that a {\em frame} is a complete lattice $L$ satisfying the
infinite distributive law $a\we\bigvee S=\bigvee\{a\we s\st s\in
S\}$, for every $a\in L$ and every $S\sbe L$.

\begin{fact}\label{frlid}{\rm (\cite{D2})}
Let $(A,\rho,\BBBB)$ be an LCA. Then the poset
$(I(A,\rho,\BBBB),\le)$ of all $\d$-ideals of $(A,\rho,\BBBB)$
(see \ref{lideal}) is a frame.
\end{fact}

\doc It is well known that the set $Idl(A)$ of all ideals of a
distributive lattice forms a frame under the inclusion ordering
(see, e.g., \cite{J}). It is easy to see that the join in
$(Idl(A),\sbe)$ of a family of $\d$-ideals  is a $\d$-ideal and
hence it is the join of this family in $(I(A,\rho,\BBBB),\le)$.
The meet in $(Idl(A),\sbe)$ of a finite family of $\d$-ideals is
also a $\d$-ideal and hence it is the meet of this family in
$(I(A,\rho,\BBBB),\le)$. Therefore, $(I(A,\rho,\BBBB),\le)$ is a
frame. Note that the meet of an infinite family of $\d$-ideals in
$(I(A,\rho,\BBBB),\le)$ is not obliged to coincide with the meet
of the same family in $(Idl(A),\sbe)$.
\sqs

\begin{theorem}\label{opensetsfr}{\rm (\cite{D2})}
Let $(A,\rho,\BBBB)$ be an LCA, $X=\Psi^a(A,\rho,\BBBB)$ and
$\OO(X)$ be the frame of all open subsets of\/ $X$. Then there
exists a frame isomorphism
$$\iota:(I(A,\rho,\BBBB),\le)\lra (\OO(X),\sbe),$$
where $(I(A,\rho,\BBBB),\le)$ is the frame of all $\d$-ideals of
$(A,\rho,\BBBB)$. The  isomorphism $\iota$ sends the set
$PI(A,\rho,\BBBB)$ of all principal $\d$-ideals of
$(A,\rho,\BBBB)$ onto the set of those regular open subsets of $X$
whose complements are in $\lag(A)$. In particular, if
$(A,\rho,\BBBB)$ is a CLCA, then $\iota(PI(A,\rho,\BBBB))=RO(X)$.
\end{theorem}

\doc  Let $I$ be a $\d$-ideal. Put $\iota(I)=\bigcup\{\l^g_A(a)\st a\in
I\}$. Then $\iota(I)$ is an open subset of $X$. Indeed, for every
$a\in I$ there exists $b\in I$ such that $a\ll b$. Then
$\l^g_A(a)\sbe\int_X(\l^g_A(b))\sbe \l^g_A(b)\sbe \iota(I)$. Hence
$\iota(I)$ is an open subset of $X$. Therefore $\iota$ is a
function from $I(A,\rho,\BBBB)$ to $\OO(X)$. Let $U\in\OO(X)$. Set
$\BBBB_U=\{b\in\BBBB\st\l^g_A(b)\sbe U\}$. Then, using
(\ref{eel}), regularity of $X$ and (III), it is easy to see that
$\BBBB_U$ is a $\d$-ideal of $(A,\rho,\BBBB)$ and
$\iota(\BBBB_U)=U$. Hence, $\iota$ is a surjection. We will show
that $\iota$ is an injection as well. Indeed, let $I_1,I_2\in
I(A,\rho,\BBBB)$ and $\iota(I_1)=\iota(I_2)$. Set $\iota(I_1)=W$
and put $\BBBB_W=\{b\in\BBBB\st\l^g_A(b)\sbe W\}$. Then,
obviously, $I_1\sbe\BBBB_W$. Further, if $b\in\BBBB_W$ then
$\l^g_A(b)\sbe W$ and $\l^g_A(b)$ is compact. Since $I_1$ is a
$\d$-ideal, $\OM=\{\int(\l^g_A(a))\st a\in I_1\}$ is an open cover
of $W$ and, hence, of $\l^g_A(b)$. Thus there exists a finite
subfamily $\{\int(\l^g_A(a_1)),\ldots,\int(\l^g_A(a_k))\}$ of
$\OM$ such that $\l^g_A(b)\sbe\bigcup\{\l^g_A(a_i)\st
i=1,\ldots,k\}=\l^g_A(\bigvee\{a_i\st i=1,\ldots,k\})$. This
implies that $b\le \bigvee\{a_i\st i=1,\ldots,k\}$ and hence $b\in
I_1$. So, we have proved that $I_1=\BBBB_W$. Analogously we can
show that $I_2=\BBBB_W$. Thus $I_1=I_2$. Therefore, $\iota$ is a
bijection. It is obvious that if $I_1,I_2\in I(A,\rho,\BBBB)$ and
$I_1\le I_2$ then $\iota(I_1)\sbe \iota(I_2)$. Conversely, if
$\iota(I_1)\sbe \iota(I_2)$ then $I_1\le I_2$. Indeed, if
$\iota(I_i)=W_i,\ i=1,2$, then, as we have already seen,
$I_i=\BBBB_{W_i},\ i=1,2$; since $W_1\sbe W_2$ implies that
$\BBBB_{W_1}\sbe\BBBB_{W_2}$, we get that $I_1\le I_2$. So,
$\iota:(I(A,\rho,\BBBB),\le)\lra(\OO(X),\sbe)$ is an isomorphism
of posets. This implies that $\iota$ is also a frame isomorphism.

Let $U$ be a regular open subset of $X$, $F=X\stm U$ and let there
exists $a\in A$ such that $F=\l^g_A(a)$. Put
$\BBBB_U=\{b\in\BBBB\st \l^g_A(b)\sbe U\}$. Then, as we have
already seen, $\BBBB_U$ is a $\d$-ideal and $\iota(\BBBB_U)=U$.
Since $F\in RC(X)$, we have that $U=X\stm
F=\int(F^*)=\int(\lag(a^*))$. Thus $\BBBB_U=\{b\in\BBBB\st
b\ll_\rho a^*\}$. Hence $\BBBB_U$ is a principal $\d$-ideal.

Conversely, if $I$ is a principal $\d$-ideal then $U=\iota(I)$ is
a regular open set in $X$ such that $X\stm U\in\lag(A)$. Indeed,
let $a\in A$ and $I=\{b\in \BBBB\st b\ll_\rho a\}$. It is enough
to prove that $X\stm U=\l^g_A(a^*)$. If $b\in I$ then $b(-\rho)
a^*$ and hence $\l^g_A(b)\cap\l^g_A(a^*)=\ems$. Thus $U\sbe X\stm
\l^g_A(a^*)$. If $\s\in X\stm\l^g_A(a^*)$ then, by (\ref{eel}),
there exists $b\in \BBBB$ such that $\s\in\l^g_A(b)\sbe X\stm
\l^g_A(a^*)$. Since, by (\ref{intha}) and (\ref{inthal}), $X\stm
\l^g_A(a^*)=\int_X(\l^g_A(a))$, we get that $b\ll_\rho a$.
Therefore $b\in I$ and hence $\s\in U$. This means that
$X\stm\l^g_A(a^*)\sbe U$.
\sqs

\begin{defi}\label{dval}{\rm (De Vries \cite{dV})}
\rm Let  $\HC$ be the category of all compact Hausdorff spaces and
all continuous maps between them.

Let $\VAL$ be the category whose objects are all complete NCAs and
whose morphisms are all functions $\p:(A,C)\lra (B,C\ap)$ between
the objects of $\VAL$ satisfying the conditions:

\smallskip

\noindent(DVAL1) $\p(0)=0$;\\
(DVAL2) $\p(a\we b)=\p(a)\we \p(b)$, for all $a,b\in A$;\\
(DVAL3) If $a, b\in A$ and $a\ll_C b$, then $(\p(a^*))^*\ll_{C\ap}
\p(b)$;\\
(DVAL4) $\p(a)=\bigvee\{\p(b)\st b\ll_{C} a\}$, for every $a\in
A$,

\medskip

{\noindent}and let the composition $``\diamond$" of two morphisms
$\p_1:(A_1,C_1)\lra (A_2,C_2)$ and $\p_2:(A_2,C_2)\lra (A_3,C_3)$
of $\VAL$ be defined by the formula
\begin{equation}\label{diamc}
\p_2\diamond\p_1 = (\p_2\circ\p_1)\cuk,
\end{equation}
 where, for every
function $\psi:(A,C)\lra (B,C\ap)$ between two objects of $\VAL$,
$\psi\cuk:(A,C)\lra (B,C\ap)$ is defined as follows:
\begin{equation}\label{cukfc}
\psi\cuk(a)=\bigvee\{\psi(b)\st b\ll_{C} a\},
\end{equation}
for every $a\in A$.
\end{defi}

De Vries \cite{dV} proved the following duality theorem:

\begin{theorem}\label{dvth}{\rm (\cite{dV})}
The categories $\HC$ and $\VAL$ are dually equivalent. In more
details, let $\Phi^t:\HC\lra\VAL$ be the contravariant functor
defined by $\Phi^t(X,\tau)=(RC(X,\tau),\rho_X)$, for every
$X\in\card\HC$, and $\Phi^t(f)(G)=\cl(f\inv(\int(G)))$, for every
$f\in\HC(X,Y)$ and every $G\in RC(Y)$, and let
$\Phi^a:\VAL\lra\HC$ be the contravariant functor defined by
$\Phi^a(A,C)=\Psi^a(A,C)$, for every $(A,C)\in\card\VAL$, and
$\Phi^a(\p)(\s\ap)=\{a\in A\st$if $b\in A$ and $b\ll_C a^*$ then
$(\p(b))^*\in\s\ap\}$, for every $\p\in\VAL((A,C),(B,C\ap))$ and
for every $\s\ap\in\CL(B,C\ap)$; then  $\l:
Id_{\,\VAL}\lra\Phi^t\circ\Phi^a$, where $\l(A,C)=\l_{(A,C)}$ (see
(\ref{h}) and (\ref{isom}) for the notation $\l_{(A,C)}$) for
every $(A,C)\in\card\VAL$, and
$t:Id_{\,\HC}\lra\Phi^a\circ\Phi^t$, where $t(X)=t_X$ (see
(\ref{homeo}) for the notation $t_X$) for every $X\in\card\HC$,
are natural isomorphisms.
\end{theorem}

In \cite{dV}, de Vries uses the regular open sets instead of
regular closed sets, as we do, so that we present here the
translations of his definitions for the case of regular closed
sets.

\begin{defi}\label{pal}
\rm We will denote by $\PLC$ the category of all locally compact
Hausdorff spaces and all perfect maps between them.

Let $\PAL$ be the category whose objects are all complete LCAs and
whose morphisms are all functions $\p:(A,\rho,\BBBB)\lra
(B,\eta,\BBBB\ap)$ between the objects of $\PAL$ satisfying the
following conditions:

\smallskip

\noindent(PAL1) $\p(0)=0$;\\
(PAL2) $\p(a\we b)=\p(a)\we \p(b)$, for all $a,b\in A$;\\
(PAL3) If $a\in\BBBB, b\in A$ and $a\llx b$, then $(\p(a^*))^*\lle
\p(b)$;\\
(PAL4) For every $b\in\BBBB\ap$ there exists $a\in\BBBB$ such that
$b\le\p(a)$;\\
\noindent(PAL5) If $a\in\BBBB$ then $\p(a)\in\BBBB\ap$;\\
(PAL6) $\p(a)=\bigvee\{\p(b)\st b\ll_{C_\rho} a\}$, for every
$a\in A$ (see (\ref{crho}) for $C_\rho$);

\medskip

{\noindent}let the composition $``\ast$" of two morphisms
$\p_1:(A_1,\rho_1,\BBBB_1)\lra (A_2,\rho_2,\BBBB_2)$ and
$\p_2:(A_2,\rho_2,\BBBB_2)\lra (A_3,\rho_3,\BBBB_3)$ of $\PAL$ be
defined by the formula
\begin{equation}\label{diam}
\p_2\ast\p_1 = (\p_2\circ\p_1)\gek,
\end{equation}
 where, for every
function $\psi:(A,\rho,\BBBB)\lra (B,\eta,\BBBB\ap)$ between two
objects of $\PAL$, $\psi\gek:(A,\rho,\BBBB)\lra (B,\eta,\BBBB\ap)$
is defined as follows:
\begin{equation}\label{cukf}
\psi\gek(a)=\bigvee\{\psi(b)\st b\ll_{C_\rho} a\},
\end{equation}
for every $a\in A$.

 By $\NAL$ we denote the full
subcategory of $\PAL$ having as objects all CNCAs (i.e., those
CLCAs $(A,\rho,\BBBB)$ for which $\BBBB=A$).
\end{defi}

Note that the categories $\VAL$ and $\NAL$ are isomorphic (it can
be even said that they are identical) because the axiom (PAL5) is
tri\-vially fulfilled in the category $\VAL$ (indeed, all elements
of its objects are bounded), the axiom (PAL4) follows immediately
from the obvious fact that $\p(1)=1$ for every $\VAL$-morphism
$\p$, and the compositions are the same.

In \cite{D}, the following generalization of de Vries Duality
Theorem was proved:

\begin{theorem}\label{gendv}{\rm (\cite{D})}
The categories $\PLC$ and $\PAL$ are dually equivalent. In more
details, let $\Psi^t:\PLC\lra\PAL$ and $\Psi^a:\PAL\lra\PLC$ be
the contravariant functors extending the Roeper correspondences
$\Psi^t$  and $\Psi^a$ (see the proof of Theorem \ref{roeperl}) to
the corresponding morphisms in the following way:
$\Psi^t(f)(G)=\cl(f\inv(\int(G)))$, for every $f\in\PLC(X,Y)$ and
every $G\in RC(Y)$, and  $\Psi^a(\p)(\s\ap)=\{a\in A\st$if $b\in
A$ and $b\ll_{C_\rho} a^*$ then $(\p(b))^*\in\s\ap\}$, for every
$\p\in\PAL((A,\rho,\BBBB),(B,\eta,\BBBB\ap))$ and for every
$\s\ap\in\Psi^a(B,\eta,\BBBB\ap)$; then  $\l^g:
Id_{\,\PAL}\lra\Psi^t\circ\Psi^a$, where
$\l^g(A,\rho,\BBBB)=\l_A^g$ for every $(A,\rho,\BBBB)\in\card\PAL$
(see (\ref{hapisomn}) for the notation $\l_A^g$), and
$t:Id_{\,\PLC}\lra\Psi^a\circ\Psi^t$, where $t(X)=t_X$ for every
$X\in\card\PLC$ (see (\ref{homeo}) for the notation $t_X$), are
natural isomorphisms.
\end{theorem}

We are now going to  generalize Theorem \ref{gendv} (and hence,
the de Vries Duality Theorem as well).

\begin{defi}\label{dhc}
\rm Let  $\HLC$ be the category of all locally compact Hausdorff
spaces and all continuous maps between them.

Let $\DLC$ be the category whose objects are all complete LCAs and
whose morphisms are all functions $\p:(A,\rho,\BBBB)\lra
(B,\eta,\BBBB\ap)$ between the objects of $\DLC$ satisfying
conditions (PAL1)-(PAL4) (which will be denoted also by (DLCi) for
$i=1\div 4$) and the following constrain:

\smallskip

\noindent(DLC5) $\p(a)=\bigvee\{\p(b)\st b\in\BBBB, b\llx a\}$,
for every $a\in A$;

\medskip

{\noindent}let the composition $``\diamond$" of two morphisms
$\p_1:(A_1,\rho_1,\BBBB_1)\lra (A_2,\rho_2,\BBBB_2)$ and
$\p_2:(A_2,\rho_2,\BBBB_2)\lra (A_3,\rho_3,\BBBB_3)$ of\/ $\DLC$
be defined by the formula
\begin{equation}\label{diamcon}
\p_2\diamond\p_1 = (\p_2\circ\p_1)\cuk,
\end{equation}
 where, for every
function $\psi:(A,\rho,\BBBB)\lra (B,\eta,\BBBB\ap)$ between two
objects of\/ $\DLC$, $\psi\cuk:(A,\rho,\BBBB)\lra
(B,\eta,\BBBB\ap)$ is defined as follows:
\begin{equation}\label{cukfcon}
\psi\cuk(a)=\bigvee\{\psi(b)\st b\in \BBBB, b\llx a\},
\end{equation}
for every $a\in A$.

 By\/ $\DHC$ we denote the full
subcategory of\/ $\DLC$ having as objects all CNCAs (i.e., those
CLCAs $(A,\rho,\BBBB)$ for which $\BBBB=A$).

(We used here the same notations as in Definition \ref{dval} for
the composition between the morphisms of the category $\DLC$ and
for the functions of the type $\psi\cuk$ because the NCAs can be
regarded as those LCAs $(A,\rho,\BBBB)$ for which $A=\BBBB$, and
hence the right sides of the formulas (\ref{cukfcon}) and
(\ref{cukfc}) coincide in the case of NCAs.)
\end{defi}

The fact that $\DLC$ is indeed a category will be proved in the
next section.

\smallskip

It is easy to show that  condition (DLC3) in Definition \ref{dhc}
can be replaced by the following one:

\medskip

\noindent(DLC3') If $a,b\in\BBBB$ and $a\llx b$, then
$(\p(a^*))^*\lle \p(b)$.

\medskip

Indeed, clearly, condition (DLC3) implies condition (DLC3').
Conversely, if $a\in\BBBB$, $b\in A$ and $a\llx b$, then, by
(BC1), there exists $c\in\BBBB$ such that $a\llx c\llx b$. Now,
using (DLC3'), we get that $(\p(a^*))^*\lle\p(c)$. Since, by
condition (DLC2), $\p$ is a monotone function, we get that
$(\p(a^*))^*\lle\p(b)$. So, conditions (DLC2) and (DLC3') imply
condition (DLC3).

Moreover, as we will see later, condition (DLC3) can be replaced
even with the following strong condition:

\medskip

\noindent(DLC3S) If $a,b\in A$ and $a\llx b$, then
$(\p(a^*))^*\lle \p(b)$.
\medskip

 We will now prove a simple lemma (almost all of which is,
 in fact, from \cite{dV} and \cite{D}), where
 some immediate consequences
  of the axioms (DLC1)-(DLC5) are listed:

  \begin{lm}\label{pf1}
Let $(A,\rho,\BBBB)$ and $(B,\eta,\BBBB\ap)$ be  CLCAs and
$\p:A\lra B$ be a function between them. Then:

\smallskip

\noindent(a) If $\p$ satisfies condition (DLC2) then $\p$ is an
order preserving function;\\
(b) If $\p$ satisfies conditions (DLC1) and (DLC2) then
$\p(a^*)\le (\p(a))^*$, for every $a\in A$; hence, if $\p$
satisfies  conditions (DLC1)-(DLC3), then for every $a\in\BBBB$
and every $b\in A$ such that $a\llx b$, $\p(a)\lle\p(b)$;\\
(c) If $\p$ satisfies conditions (DLC2) and (DLC4) (or (DLC1) and
(DLC3)) then $\p(1_A)=1_B$;\\
(d) If $\p$ satisfies condition (DLC2) then $\p\cuk$ satisfies
conditions (DLC2) and (DLC5) (see (\ref{cukf}) for $\p\cuk$);\\
(e) If $\p$ satisfies condition (DLC5) then $\p=\p\cuk$;\\
(f) If $\p$ satisfies condition (DLC2) then
$(\p\cuk)\cuk=\p\cuk$;\\
(g) If $\p$ is a monotone function then, for every $a\in A$,
$\p\cuk(a)\le\p(a)$.
%
%(h) If $\p$ is an LCA-isomorphism then
%$(\p\circ\a\circ\b)\cuk=\p\circ(\a\circ\b)\cuk$, where $\a$ and
%$\b$ are functions between CLCAs for which the composition
%$\p\circ\a\circ\b$ exists.
\end{lm}

\doc The properties (a), (b), (e) and (g) are clearly fulfilled, and (f)
follows from (d) and (e).

\medskip

\noindent(c) Using consecutively (a) and (DLC4), we get that
$\p(1_A)\ge \bigvee\{\p(b)\st b\in\BBBB\}\ge\bigvee\{b\ap\st
b\ap\in\BBBB\ap\}=1_B$. Thus, $\p(1_A)=1_B$. The assertion in
brackets can be obtained easily applying (DLC3) and (DLC1) to the
inequality $0_A\ll_\rho 0_A$.

\medskip

\noindent(d) By (a), for every $a\in A$, $\p\cuk(a)\le\p(a)$. Let
$a\in A$. If $c\in\BBBB$ and $c\llx a$ then there exists $d_c\in
\BBBB$ such that $c\llx d_c\llx a$; hence $\p(c)\le\p\cuk(d_c)$.
Now, $\p\cuk(a)=\bigvee\{\p(c)\st c\in\BBBB, c\llx
a\}\le\bigvee\{\p\cuk(d_c)\st c\in\BBBB,  c\llx
a\}\le\bigvee\{\p\cuk(e)\st e\in\BBBB,  e\llx
a\}\le\bigvee\{\p(e)\st e\in\BBBB,  e\llx a\}=\p\cuk(a)$. Thus,
$\p\cuk(a)=\bigvee\{\p\cuk(e)\st e\in\BBBB, e\llx a\}$. So,
$\p\cuk$ satisfies (DLC5). Further, let $a,b\in A$. Then
$\p\cuk(a)\we\p\cuk(b)=\bigvee\{\p(d)\we\p(e)\st d,e\in\BBBB,
d\llx a, e\llx b\}=\bigvee\{\p(d\we e)\st d,e\in\BBBB,  d\llx a,
e\llx b\}=\bigvee\{\p(c)\st c\in\BBBB,  c\llx a\we b\}=\p\cuk(a\we
b)$. So, (DLC2) is fulfilled. \sqs

\medskip

Obviously, the assertions (a), (b) and (c) of the above lemma
remain true also in the case when $(A,\rho,\BBBB)$ and
$(B,\eta,\BBBB\ap)$ are LCAs.

As we shall prove in the next section,  condition (DLC3) in
Definition \ref{dhc} can be also replaced by any of the following
two constrains:

\medskip

\noindent(LC3) If, for $i=1,2$, $a_i\in\BBBB$, $b_i\in A$  and
$a_i\llx b_i$, then $\p(a_1\vee a_2)\lle \p(b_1)\vee \p(b_2)$.

\medskip

\noindent(LC3S) If, for $i=1,2$, $a_i,b_i\in A$  and $a_i\llx
b_i$, then $\p(a_1\vee a_2)\lle \p(b_1)\vee \p(b_2)$.

\medskip

\noindent  Now, we will only  show  that the following proposition
holds:

\begin{pro}\label{axlc3}
Conditions (DLC1)-(DLC3) in Definition \ref{dhc} imply condition
(LC3).
\end{pro}

\doc Let, for $i=1,2$, $a_i\in\BBBB$,  $b_i\in A$  and $a_i\llx
b_i$. Then, by (DLC3), $(\p(a_i^*))^*\lle\p(b_i)$, where $i=1,2$.
Hence $(\p(b_i))^*\lle\p(a_i^*)$, for $i=1,2$. Thus, using
consecutively (DLC2) and Lemma \ref{pf1}(b), we get that
$(\p(b_1))^*\we(\p(b_2))^*\lle\p(a_1^*)\we\p(a_2^*)=\p(a_1^*\we
a_2^*)=\p((a_1\vee a_2)^*)\le(\p(a_1\vee a_2))^*$. Therefore,
$\p(a_1\vee a_2)\lle \p(b_1)\vee \p(b_2)$. \sqs

\begin{rem}\label{remlc3}
\rm Obviously, condition (LC3) implies the second assertion of
Lemma \ref{pf1}(b).
\end{rem}

\begin{theorem}\label{lccont}{\rm (Main Theorem)}
The categories $\HLC$ and\/ $\DLC$ are dually equivalent. In more
details, let
$$\LAM^t:\HLC\lra\DLC \mbox{ and }\LAM^a:\DLC\lra\HLC$$
be the contravariant functors extending, respectively,  the Roeper
correspondences $\Psi^t:\card{\HLC}\lra\card{\DLC}$  and\/
$\Psi^a:\card{\DLC}\lra\card{\HLC}$ (see (\ref{psit1}),
(\ref{phiapcn}) and (\ref{phiapc}) for the notations) to the
corresponding morphisms of the categories $\HLC$ and\/ $\DLC$ in
the following way:
$$\LAM^t(f)(G)=\cl(f\inv(\int(G))),$$
for every $f\in\HLC(X,Y)$ and every $G\in RC(Y)$, and
$$\LAM^a(\p)(\s\ap)\cap\BBBB=\{a\in \BBBB\st \mbox{if } b\in A
\mbox{ and } a\llx b \mbox{ then }\p(b)\in\s\ap\},$$
for every $\p\in\DLC((A,\rho,\BBBB),(B,\eta,\BBBB\ap))$ and for
every $\s\ap\in\LAM^a(B,\eta,\BBBB\ap)$; then
$$\l^g: Id_{\,\DLC}\lra\LAM^t\circ\LAM^a, \mbox{ where }
\l^g(A,\rho,\BBBB)=\l_A^g$$
for every $(A,\rho,\BBBB)\in\card\DLC$ (see (\ref{hapisomn}) for
the notation $\l_A^g$), and
%%%
$$t^l:Id_{\,\HLC}\lra\LAM^a\circ\LAM^t, \mbox{ where }t^l(X)=t_X$$
for every $X\in\card\HLC$ (see (\ref{homeo}) for the notation
$t_X$), are natural isomorphisms.
\end{theorem}

The proof of this theorem will be presented in the next section.

\smallskip

It is clear that the categories $\DHC$ and $\VAL$ (see \ref{dhc}
and \ref{dval} for their definitions) are isomorphic (it can be
even said that they are identical). Hence, using the fact that
compact spaces correspond to CNCAs and conversely (see Theorem
\ref{roeperl}), we get, by Theorem \ref{lccont}, that the
ca\-tegories $\VAL$ and $\HC$ are dually equivalent. Moreover, the
definitions of the corresponding duality functors coincide.
Indeed, it is obvious that the definitions of the contravariant
functor $\Phi^t$ and the restriction of the contravariant functor
$\LAM^t$ to the subcategory $\HC$ of the category $\HLC$ coincide.
Further, we need to show that the contravariant functor $\Phi^a$
and the restriction of the contravariant functor $\LAM^a$ to the
subcategory $\DHC$ of the category $\DLC$ coincide. Let
$\p\in\DHC((A,\rho,A),(B,\eta,B))$($=\VAL((A,\rho),(B,\eta))$) and
$\s\ap\in\Psi^a(B,\eta,B)$.  Then set
$\s_{\Phi}=\Phi^a(\p)(\s\ap)=\{a\in A \st (\fa b\in A)[(b\ll
a^*)\rightarrow ((\p(b))^*\in\s\ap)]\}$. Obviously,
$\s_{\Phi}=\{a\in A \st (\fa b\in A)[(a\ll b)\rightarrow
((\p(b^*))^*\in\s\ap)]\}$. Set $\s_{\LAM}=\LAM^a(\p)(\s\ap)=\{a\in
A \st (\fa b\in A)[(a\ll b)\rightarrow (\p(b)\in\s\ap)]\}$. Then
$\s_{\LAM}\sbe\s_{\Phi}$. Indeed, by Lemma \ref{pf1}(b),
$\p(b^*)\le (\p(b))^*$; hence $\p(b)\le (\p(b^*))^*$; therefore,
if $a\in\s_{\LAM}$ then $a\in\s_{\Phi}$. Now, Fact \ref{fact29}
implies that $\s_{\LAM}=\s_{\Phi}$.
  Thus, de Vries Duality Theorem is a corollary
of Theorem \ref{lccont}.

\begin{defi}\label{defplc}
\rm Let $\PDLC$ be the subcategory of the category $\DLC$ whose
objects are all CLCAs and whose morphisms are all $\DLC$-morphisms
satisfying condition (PAL5).
\end{defi}

We will now show that the categories $\PDLC$ and $\PAL$ are
isomorphic (it can be even said that they are identical).

\begin{lm}\label{cuklm}
Let $(A,\rho,\BBBB)$ and $(B,\eta,\BBBB\ap)$ be  CLCAs and
 $\psi:A\lra B$ satisfy conditions (DLC2), (DLC4). Then, for
 every $a\in A$, $\psi\cuk(a)=\psi\gek(a)$. (See (\ref{cukfcon}) and
 (\ref{cukf}) for the notations  $\psi\cuk$ and $\psi\gek$.)
\end{lm}

\doc   Obviously, if $b\in\BBBB$ and $b\llx a$ then
$b\llcr a$. Thus, $\psi\cuk(a)\le \psi\gek(a)$.

Let now $b\in A$ and $b\llcr a$. Then $b\llx a$. We have, by
(\ref{aunionp}), that $\psi(b)=\bigvee\{c\,\ap\in\BBBB\ap\st
c\,\ap\lle\psi(b)\}$. Let $c\,\ap\in\BBBB\ap$ and
$c\,\ap\lle\psi(b)$. Then, by (DLC4), there exists $c\in\BBBB$
such that $c\,\ap\le\psi(c)$. Now, (DLC2) implies that
$c\,\ap\le\psi(b\we c)$. Set $d=b\we c$. Then $d\in\BBBB$, $d\le
b\llx a$ (and, hence, $d\llx a$), $c\,\ap\le\psi(d)\le
\psi\cuk(a)$. Thus, $\psi(b)\le \psi\cuk(a)$. We conclude that
$\psi\gek(a)\le \psi\cuk(a)$. So, $\psi\cuk(a)=\psi\gek(a)$. \sqs

\begin{rem}\label{cuklmrem}
\rm The proof of Lemma \ref{cuklm} shows even that
$\psi\cuk(a)=\psi\gek(a)=\bigvee\{\psi(b)\st b\llx a\}$.
\end{rem}

\begin{cor}\label{axpal}
The system of axioms (PAL1)-(PAL6) is equivalent
 to the system of axioms (PAL1)-(PAL5), (DLC5).
 \end{cor}

 \begin{cor}\label{comppal}
The compositions in\/ $\PAL$  and\/ $\PDLC$ coincide.
 \end{cor}

 \doc Let $\p_1,\p_2$ be $\PAL$-morphisms. We have to show that
 $\p_2\ast\p_1=\p_2\diamond\p_1$. Set $\psi=\p_2\circ\p_1$. Then,
 using \ref{pf1}(a), one obtains immediately that $\psi$
 satisfies conditions (DLC2), (DLC4). Now we can apply \ref{cuklm}.
 \sqs

All this proves the following assertion:

\begin{pro}\label{corplc}
The categories\/ $\PAL$ and\/ $\PDLC$ are isomorphic.
\end{pro}

We will now show, using our Main Theorem, that the categories
$\PDLC$ and $\PLC$ are dually equivalent. In this way, we will
obtain Theorem \ref{gendv}  as a corollary of Theorem
\ref{lccont}.

\begin{theorem}\label{lcper}
The categories\/ $\PLC$ and\/ $\PDLC$ are dually equivalent.
\end{theorem}

\doc  We will show that the restriction to the subcategory $\PDLC$
of the category $\DLC$ of the duality
functor $\LAM^a:\DLC\lra\HLC$, constructed in Theorem
\ref{lccont},  is the desired duality functor between the
categories $\PDLC$ and $\PLC$.

Let $f\in\PLC(X,Y)$. Then, by Theorem \ref{lccont},
$\p_f=\LAM^t(f)$ is an element of the set
$\DLC((RC(Y),\rho_Y,CR(Y)), (RC(X),\rho_X,CR(X)))$.
%%%Set $\p=\p_f$.
We will show that $\p_f$ is a $\PDLC$-morphism, i.e. that $\p_f$
satisfies, in addition, condition (PAL5). So, let $G\in CR(Y)$.
Then $\p_f(G)=\cl_X(f\inv(\int(G)))$ (see Theorem \ref{lccont}),
and hence $\p_f(G)\sbe f\inv(G)$. Since $f$ is a perfect map,
$f\inv(G)$ is a compact subset of $X$. Thus $\p_f(G)\in CR(X)$.
Therefore, condition (PAL5) is fulfilled.

Let now $\p\in\PDLC((A,\rho,\BBBB),(B,\eta,\BBBB\ap))$ and
$f_\p=\LAM^a(\p)$ (see Theorem \ref{lccont}). Let
$X=\Psi^a(A,\rho,\BBBB)$ and $Y=\Psi^a(B,\eta,\BBBB\ap)$. Then, by
Theorem \ref{lccont}, $f_\p:Y\lra X$ is a continuous map. We will
show that $f_\p$ is a perfect map. Let us first prove that for
every $a\in\BBBB$, $f_\p\inv(\lag(a))$ is a compact subset of $Y$.
Indeed, let $a\in\BBBB$. Then, by condition (BC1) (see
\ref{locono}), there exists $b\in\BBBB$ such that $a\llx b$. We
will show that $f_\p\inv(\lag(a))\sbe\lbg(\p(b))$.  So, let
$\s\ap\in f_\p\inv(\lag(a))$. Then $f_\p(\s\ap)\in\lag(a)$. Hence
$a\in\BBBB\cap f_\p(\s\ap)$. This implies, by Theorem
\ref{lccont}, that $\p(b)\in\s\ap$. Thus $\s\ap\in\lbg(\p(b))$.
Since, by (PAL5), $\p(b)\in\BBBB\ap$, we have that $\lbg(\p(b))$
is compact. Hence $f_\p\inv(\lag(a))$ is compact.
 So, for every $a\in\BBBB$, $f_\p\inv(\lag(a))$ is a compact subset of $Y$.
Now, using the fact that the family $\{\int(\lag(a))\st a\in
\BBBB\}$ is an open base of $X$, we conclude that $f_\p\inv(K)$ is
compact for every compact subset $K$ of $X$. Then the local
compactness of $X$ implies that $f_\p$ is a perfect map.
Therefore, we have proved that $f_\p\in\PLC(Y,X)$.

The rest follows from Theorem \ref{lccont}. \sqs

\begin{nota}\label{connlc}
\rm We will denote by $\CHLC$ the full subcategory of the
ca\-tegory $\HLC$ whose objects are all connected locally compact
Hausdorff spaces, and by $\CDLC$ the full subcategory of the
category $\DLC$ whose objects are all connected CLCAs.
\end{nota}

\begin{theorem}\label{lcconn}
The categories\/ $\CHLC$ and\/ $\CDLC$ are dually equivalent.
\end{theorem}

\doc It follows immediately from Theorem \ref{lccont} and Fact
\ref{confact}. \sqs

Finally, we will show that the categories $\SKLC$ and $\SKAL$,
introduced in the paper \cite{D1} in which the Fedorchuk Duality
Theorem \cite{F} is generalized, are (non full) subcategories of
the categories, respectively, $\HLC$ and $\DLC$ and the
restrictions of the duality functors $\LAM^a$ and $\LAM^t$ to them
coincide with the corresponding duality functors described in
\cite{D1}.

Let us first formulate  Fedorchuk Duality Theorem. We will need
the following definition:  a continuous map $f:X\lra Y$ is called
{\em quasi-open\/} (\cite{MP}) if for every non-empty open subset
$U$ of $X$, $\int(f(U))\nes$ holds.

We will denote by $\SHC$ the category of all compact Hausdorff
spaces and all quasi-open maps between them.

Let $\SAC$ be the category whose objects are all complete normal
contact algebras and whose morphisms $\p:(A,C)\lra (B,C\ap)$ are
all complete Boolean homomorphisms $\p:A\lra B$ satisfying the
following condition:

\smallskip

\noindent (F1) For all $a,b\in A$,  $\p(a)C\ap\p(b)$ implies
$aCb$.

\begin{theorem}\label{dcompn}{\rm (Fedorchuk \cite{F})}
The categories $\SHC$ and $\SAC$ are dually equi\-valent.
\end{theorem}

We are now going to formulate a generalization (presented in
\cite{D1}) of the last theorem.

 Let $\p:A\lra B$ be an order-preserving map between posets.
Then a map $\p_\LAM :B\lra A$ is called a {\em left adjoint} to
$\p$ if it is the unique order-preserving map such that, for all
$a\in A$ and all $b\in B$, $b\le \p(a)$ iff $\p_\LAM (b)\le a$
(i.e. the pair $(\pl,\p)$ forms a Galois connection between posets
$B$ and $A$).

 Recall
that a function $f:X\lra Y$ is called {\em skeletal}\/ (\cite{MR})
if
\begin{equation}\label{ske}
\int(f\inv(\cl (V)))\sbe\cl(f\inv(V))
\end{equation}
for every open subset  $V$  of $Y$.

Let $\SKLC$ be the category of all locally compact Hausdorff
spaces and all continuous skeletal maps between them.

Let $\SKAL$ be the category whose objects are all complete local
contact algebras
%(see \ref{locono})
 and whose morphisms   $\p:(A,\rho,\BBBB)\lra
(B,\eta,\BBBB\ap)$   are all
  complete Boolean homomorphisms $\p:A\lra B$  satisfying the
following conditions:

\smallskip

\noindent (L1) $\fa a,b\in A$, $\p(a)\eta\p(b)$ implies $a\rho
b$;\\
(L2) $b\in\BBBB\ap$ implies $\pl(b)\in\BBBB$ (where $\pl$ is the
left adjoint to $\p$).

\smallskip

Let us  note  that (L1) is equivalent to the following condition:

\smallskip

\noindent (EL1) $\fa a,b\in B$, $a\eta b$ implies
$\pl(a)\rho\pl(b)$.

\begin{theorem}\label{maintheoremnk}{\rm  (\cite{D1})}
The categories $\SKLC$ and $\SKAL$ are dually equi\-valent. In
more details, let
$$\Psi_1^t:\SKLC\lra\SKAL \mbox{ and }\Psi_1^a:\SKAL\lra\SKLC$$
 be
two contravariant functors extending the Roeper correspondences
$\Psi^t$  and $\Psi^a$ (see the proof of Theorem \ref{roeperl}) to
the corresponding morphisms in the following way:
$\Psi_1^t(f)(G)=\cl(f\inv(\int(G)))$, for every $f\in\SKLC(X,Y)$
and every $G\in RC(Y)$, and
 for every
$\p\in\SKAL((A,\rho,\BBBB),(B,\eta,\BBBB\ap))$ and for every
bounded ultrafilter  $u$ in $B$  (see \ref{boundcl})
\begin{equation}\label{psi1a}
\Psi_1^a(\p)(\s_u)=\s_{\p\inv(u)},
\end{equation}
where $\s_{\p\inv(u)}$ is a cluster in $(A,C_\rho)$ (see
\ref{Alexprn}, (\ref{sigmau}) and \ref{uniqult} for  $C_\rho$ and
$\s_u$, and note that, by \ref{conclustth}, any bounded cluster
$\s$ in $(B,\eta,\BBBB\ap)$ can be written in the form $\s_u$ for
some bounded ultrafilter $u$ in $B$); then $\l^g:
Id_{\,\SKAL}\lra\Psi_1^t\circ\Psi_1^a$, where
$\l^g(A,\rho,\BBBB)=\l_A^g$, for every
$(A,\rho,\BBBB)\in\card\SKAL$ (see (\ref{hapisomn}) for the
notation $\l_A^g$), and $t:Id_{\,\SKLC}\lra\Psi_1^a\circ\Psi_1^t$,
where $t(X)=t_X$, for every $X\in\card\SKLC$ (see (\ref{homeo})
for the notation $t_X$), are natural isomorphisms.
\end{theorem}

Now we will prove the following proposition.

\begin{pro}\label{mainfed}
The categories  $\SKLC$ and $\SKAL$ are (non full) subcategories
of, respectively, $\HLC$ and $\DLC$. The restriction of the
contravariant functor $\LAM^a$ (respectively, $\LAM^t$) to  the
subcategory $\SKAL$ (resp., $\SKLC$) coincides with the
contravariant functor $\Psi_1^a$ (resp.,  $\Psi_1^t$).
\end{pro}

\doc  Obviously, the category $\SKLC$ is a subcategory of the
category $\HLC$ and the restriction of the contravariant functor
 $\LAM^t$ to  the subcategory $\SKLC$ coincides with the contravariant functor
 $\Psi_1^t$.

 Let $\p\in\SKAL((A,\rho,\BBBB),(B,\eta,\BBBB\ap))$. Then,
 clearly, $\p$ satisfies conditions (DLC1) and (DLC2). Let
 $a,b\in A$ and $a\llx b$. Then $a(-\rho)b^*$. Hence, by (L1),
 $\p(a)(-\eta)\p(b^*)$. Since $\p$ is a Boolean homomorphism, we
 have that $\p(b^*)=(\p(b))^*$ and $(\p(a^*))^*=\p(a)$. Thus,
 $(\p(a^*))^*\lle\p(b)$. Therefore, condition (DLC3) is satisfied.
 Let now $b\in\BBBB\ap$. Then, by (L2), $\pl(b)\in\BBBB$. Set
 $a=\pl(b)$. Since $\pl$ is a left adjoint to $\p$, we get that
 $b\le\p(a)$. So, condition (DLC4) is checked. Finally, let $a\in
 A$. Then $a=\bigvee\{b\in\BBBB\st b\llx a\}$. Since $\p$ is a
 complete Boolean  homomorphism, we conclude that
 $\p(a)=\bigvee\{\p(b)\st b\in\BBBB, b\llx a\}$.
 Thus, condition (DLC5) is satisfied. So, every $\SKAL$-morphism
 is a $\DLC$-morphism. Since the composition $\p_2\circ\p_1$ of two Boolean
 homomorphisms is a complete Boolean homomorphism, Lemma
 \ref{pf1}(e) implies that $(\p_2\circ\p_1)\cuk=\p_2\circ\p_1$.
 Hence, $\p_2\diamond\p_1=\p_2\circ\p_1$. Therefore, the category
   $\SKAL$ is a subcategory of the
category $\DLC$.

Let $\p\in\SKAL((A,\rho,\BBBB),(B,\eta,\BBBB\ap))$ and $u$ be a
bounded ultrafilter  in $B$. Then, as it is proved in \cite{D1},
$\s_{\p\inv(u)}$
 is a bounded cluster in $(A,\rho,\BBBB)$. We have that
$\Psi_1^a(\p)(\s_u)=\s_{\p\inv(u)}$ and
$\LAM^a(\p)(\s_u)\cap\BBBB=\{a\in \BBBB\st \mbox{if } b\in A
\mbox{ and } a\llx b \mbox{ then }\p(b)\in\s_u\}$. According to
Corollary \ref{bbcl} (see below),
$\Psi_1^a(\p)(\s_u)=\LAM^a(\p)(\s_u)$ iff
$\BBBB\cap\Psi_1^a(\p)(\s_u)=\BBBB\cap\LAM^a(\p)(\s_u)$. Thus, we
have to show that
\begin{equation}\label{mfe}
\BBBB\cap\s_{\p\inv(u)}=\{a\in \BBBB\st \mbox{ if } b\in A \mbox{
and } a\llx b \mbox{ then }\p(b)\in\s_u\}.
\end{equation}

So, let $a\in\BBBB\cap\s_{\p\inv(u)}$. Suppose that there exists
$b\in A$ such that $a\llx b$ and $\p(b)\nin\s_u$. Then $\p(b)\nin
u$. Hence $(\p(b))^*\in u$, i.e. $\p(b^*)\in u$. Thus
$b^*\in\p\inv(u)$. Since $a(-\rho)b^*$, we get a contradiction.

Conversely, let $a\in \BBBB$ and for all $b\in A$ such that $a\llx
b$, we have that $\p(b)\in\s_u$. We have to prove that
$a\in\s_{\p\inv(u)}$, i.e. that $a\rho b$ for all $b\in\p\inv(u)$.
Suppose that there exists $b_0\in\p\inv(u)$ such that
$a(-\rho)b_0$. Then $\p(b_0)\in u$ and $a\llx b_0^*$. By (BC1),
there exists $a_1\in\BBBB$ such that $a\llx a_1\llx b_0^*$. Hence
$a_1(-\rho)b_0$ and $b_0\le a_1^*$. Then $\p(b_0)\le\p(a_1^*)$ and
thus $\p(a_1^*)\in u$. Since $\p(a_1^*)=(\p(a_1))^*$, we get that
$\p(a_1)\nin u$. We have that $\p(a_1)\in\s_u$ (because $a\llx
a_1$).  Let $c_0\in u\cap\BBBB\ap$. Then $c_0\we\p(b_0)\in
u\cap\BBBB\ap$. Thus $\p(a_1)\eta(c_0\we\p(b_0))$. Therefore
$\p(a_1)\eta\p(b_0)$. Then, by (EL1),
$\pl(\p(a_1))\rho\pl(\p(b_0))$. Since, for every $c\in A$,
$\pl(\p(c))\le c$, we obtain that $a_1\rho b_0$, a contradiction.
So, the equality (\ref{mfe}) is established. This completes the
proof. \sqs

%--------------------------------------------------------------------
%----------------------------------------------------------------------
%%%-----------------------------------------------------------------

\section{The Proof of the Main Theorem}

 The proof of Theorem \ref{lccont}, which is our main
theorem, will be divided in several lemmas, propositions, facts
and remarks. The plan of it is the following: we begin with some
preparatory assertions; after that we show that $\DLC$ is indeed a
category; the crucial step in the proof of this fact is to show
that any function between CLCAs, which satisfies conditions
(DLC1)-(DLC5), satisfies condition (DLC3S) as well; this statement
is obtained as a corollary of
 some other assertions which are also used  later on in the last
 portion of the proof where the
 construction of the desired duality between the categories
$\DLC$ and $\HLC$ is presented.

We  start with some simple, but important for our proof,
propositions about (bounded) clusters in LCAs.

\begin{pro}\label{bbclf}{\rm (\cite{D1})}
Let $(A,\rho,\BBBB)$ be an LCA.  If $u$ is an ultrafilter in $A$
and $\s_u\cap\BBBB\nes$, then $u\cap\BBBB\nes$.
\end{pro}

\doc  By Fact \ref{bstar}, there exists $a\in\BBBB$ such that
$a^*\nin\s_u$. Then $a\in u\cap\BBBB$. \sqs

\begin{pro}\label{bbcl1}
Let $(A,\rho,\BBBB)$ be an LCA and $\s$ be a bounded cluster in
it. Then:

\smallskip

\noindent(a) If $a\in\s$ then there exists $c\in\BBBB\cap\s$ such
that $c\le a$;

\smallskip

\noindent(b) $(\fa a\in A)[(a\nin\s)\leftrightarrow((\ex
b\in\BBBB\cap\s)(b\llx a^*))]$;

\smallskip

\noindent(c) $\s=\{a\in A\st a\rho(\s\cap\BBBB)\}$.
\end{pro}

\doc (a) By Theorem \ref{conclustth}, there exists an ultrafilter $u$
in $A$ such that $a\in u$ and $\s=\s_u$. Then $u\sbe \s$. Since
$\s\cap\BBBB\nes$, Proposition \ref{bbclf} implies that there
exists $a_1\in u\cap\BBBB$. Set $c=a\we a_1$. Then $c\in
u\cap\BBBB\sbe\s\cap\BBBB$ and $c\le a$.

\smallskip

\noindent(b) Let $a\nin\s$. Then, by \ref{cluendcor}, there exists
$c\in A$ such that $c\nin\s$ and $a\llcr c$. Then $c^*\llx a^*$
and $c^*\in\s$. Hence, by (a), there exists $b\in\BBBB\cap\s$ such
that $b\llx a^*$. Conversely, if there exists $b\in\BBBB\cap\s$
such that $b\llx a^*$, then $b\llcr a^*$. Therefore, $a\nin\s$.

\smallskip

\noindent(c) This is just another form of (b).
 \sqs

 \begin{cor}\label{bbcl}
Let $(A,\rho,\BBBB)$ be an LCA and $\s_1$, $\s_2$ be two clusters
in $(A,\rho,\BBBB)$ such that\/ $\BBBB\cap\s_1=\BBBB\cap\s_2$.
Then $\s_1=\s_2$.
\end{cor}

\doc  By \ref{neogrn}, $\s_\infty=A\stm\BBBB$ is a cluster in
$(A,C_\rho)$. Hence, if  $\BBBB\cap\s_1=\BBBB\cap\s_2=\ems$, then
$\s_i\sbe\s_\infty$, for $i=1,2$. Now, \ref{fact29} implies that
$\s_1=\s_\infty=\s_2$.

Let $\BBBB\cap\s_1\nes$. Then our assertion follows from
Proposition \ref{bbcl1}(c). \sqs

Recall that if $A$ is a lattice then an element $p\in A\stm\{1\}$
is called a {\em prime element of} $A$ if for each $a,b\in A$,
$a\we b\le p$ implies that $a\le p$ or $b\le p$. We will now show
that if $(A,\rho,\BBBB)$ is a CLCA then the prime elements of
$I(A,\rho,\BBBB)$ are in a bijective correspondence with the
bounded clusters in $(A,\rho,\BBBB)$. The existence of such a
bijection follows immediately from Roeper's Theorem \ref{roeperl},
our Theorem \ref{opensetsfr} and localic duality (see, e.g.,
\cite{J}). We will present here an explicit formula for this
bijection which will be very useful later on.

\begin{pro}\label{primdid}
Let $\s$ be a bounded cluster in an LCA $(A,\rho,\BBBB)$. Then
$I=\BBBB\stm\s$ is a prime element of the frame $I(A,\rho,\BBBB)$
(see \ref{lideal} for this notation).
\end{pro}

\doc We have that $I\neq\BBBB$ because $\s\cap\BBBB\nes$.
Since $\s$ is an upper set, we get that $I$ is a lower set. Let
$a,b\in I$. Suppose that $a\vee b\in\s$. Then $a\in\s$ or
$b\in\s$, a contradiction. Hence, $a\vee b\in I$. So, $I$ is an
ideal. Let $a\in I$. Then $a\nin\s$. By \ref{cluendcor}, there
exists $c\in A$ such that $c\nin\s$ and $a\llcr c$. Thus $a\llx
c$. By (BC1), there exists $b\in\BBBB$ such that $a\llx b\llx c$.
Then $b\nin\s$, i.e. $b\in I$ and $a\llx b$. So, $I$ is a
$\d$-ideal. Let $J_1,J_2\in I(A,\rho,\BBBB)$ and $J_1\cap J_2\sbe
I$.  Suppose that, for $i=1,2$, there exists $a_i\in J_i\stm I$.
Since, for $i=1,2$, $J_i$ is a $\d$-ideal, there exists $b_i\in
J_i$ such that $a_i\llx b_i$, and hence $a_i\llcr b_i$; thus
$b_i^*\nin\s$. Then $b_1^*\vee b_2^*\nin\s$. Therefore, $b_1\we
b_2\in\s\cap I$. Since $\s\cap I=\ems$, we get a contradiction.
All this shows that $I$ is a prime element of the frame
$I(A,\rho,\BBBB)$. \sqs

\begin{pro}\label{primdfil}
 Let $(A,\rho,\BBBB)$ be an
LCA and $I$ be a prime element of the frame $I(A,\rho,\BBBB)$.
 Then
the set $V=\{a\in\BBBB\st (\ex b\in\BBBB\stm I)(b\llx a)\}$ is a
filter in $\BBBB$.
\end{pro}

\doc Set $S=\BBBB\stm I$. Then $S$ is a non-void upper set in $\BBBB$.
 Thus, $V\sbe
S$. If $a\in S$ then, by (BC1), there exists $b\in\BBBB$ such that
$a\llx b$. Then $b\in V$, i.e. $V\nes$. Obviously, $0\nin V$ and
$V$ is an upper set in $\BBBB$. Let $a,b\in V$ and suppose that
$a\we b\nin V$. Then, for every $c\in S$, $c\not\llx a\we b$.
Hence $I_a\cap I_b\sbe I$ (see \ref{dideal} for the notations).
Thus, $I_a\sbe I$ or $I_b\sbe I$. Let, e.g., $I_a\sbe I$. Since
$a\in V$, there exists $c\in S$ such that $c\llx a$. Then $c\in
I_a\cap S\sbe I\cap S=\ems$, a contradiction. Therefore, $a\we
b\in V$. So, $V$ is a filter in $\BBBB$. \sqs

\begin{pro}\label{primdidclu}
 Let $(A,\rho,\BBBB)$ be an
LCA and $I$ be a prime element of the frame $I(A,\rho,\BBBB)$.
Then there exists a unique cluster $\s$ in $(A,\rho,\BBBB)$ such
that $\s\cap\BBBB=\BBBB\stm I$; moreover, $\s=\{a\in A\st
a\rho(\BBBB\stm I)\}$. (In this case we will say that $\s$ is {\em
generated by} $I$.)
\end{pro}

\doc   By Proposition \ref{primdfil}, the set $V$ defined there is a filter in
$\BBBB$. Hence, $V\nes$ and $V$ is a filter-base in $A$. Let $F$
be the filter in $A$ generated by the filter-base $V$. Then
$F\cap\BBBB=V$ and hence $F\cap I=\ems$. Now, the famous Stone
Separation Theorem (see, e.g., \cite{J}) implies that there exists
an ultrafilter $u$ in $A$ such that $F\sbe u$ and $u\cap I=\ems$.
Set $\s=\s_u$ (see Theorem \ref{uniqult} for the notation $\s_u$).
Then $\s$ is a cluster in $(A,\rho,\BBBB)$ (i.e., $\s$ is a
cluster in the NCA $(A,C_\rho)$) and $\s=\{a\in A\st aC_\rho b$
for every  $b\in u\}$. Since $F\sbe u\sbe\s$, we have that
$V\sbe\s\cap\BBBB$. Let us show that $\s\cap I=\ems$. Indeed,
suppose that $a\in\s\cap I$. Since $I$ is a $\d$-ideal, there
exists $b\in I$ such that $a\llx b$. Then, obviously, $b\in\s\cap
I$. We have that $a(-\rho)b^*$ and $a\in\BBBB$. Hence
$a(-C_\rho)b^*$. Since $u\sbe\s$ and $a\in\s$, we get that
$b^*\nin u$. Thus $b\in u$, i.e. $b\in u\cap I$, a contradiction.
So, $\s\cap I=\ems$, i.e. $\s\cap\BBBB\sbe S=\BBBB\stm I$. We will
now prove that $\s\cap\BBBB=S$. Indeed, suppose that there exists
$a\in S\stm\s$. Then there exists $b\in u$ such that $a(-\rho)b$.
Hence $a\llx b^*$. Now, (BC1) implies that there exists
$c\in\BBBB$ such that $a\llx c\llx b^*$. Then $c\in V$ and $c\le
b^*$. Thus $c\we b=0$, which means that $c\nin u$. Since $V\sbe u$
and $c\in V$, we get a contradiction. So, $\s\cap\BBBB=S$.
Finally, Corollary \ref{bbcl} implies the uniqueness of $\s$ and
the formula  $\s=\{a\in A\st a\rho(\BBBB\stm I)\}$ follows from
\ref{bbcl1}(c). \sqs

\begin{cor}\label{bijecpecl}
For any LCA $(A,\rho,\BBBB)$, there exists a bijective
correspondence between the bounded clusters in $(A,\rho,\BBBB)$
and the prime elements of the frame $I(A,\rho,\BBBB)$ (see
\ref{lideal} for this notation).
\end{cor}

\doc It follows from \ref{primdid} and \ref{primdidclu}.
\sqs

\begin{rem}\label{bcpe}
\rm If $(A,\rho,\BBBB)$ is an LCA, then it is easy to see that
every prime ideal $J$ of $\BBBB$ (i.e. $J$ is an ideal, $J\not
=\BBBB$ and $(\fa a,b\in\BBBB)[(a\we b\in J)\rightarrow (a\in J
\mbox { or } b\in J)]$) which is a $\d$-ideal (shortly, prime
$\d$-ideal) is a prime element of the frame $I(A,\rho,\BBBB)$.
However, in contrast to the case of ideals of a lattice (where the
prime elements of the frame of all ideals of this lattice are
precisely the prime ideals of the lattice), the prime elements of
the frame $I(A,\rho,\BBBB)$ need not be prime $\d$-ideals of
$\BBBB$. Indeed, let $I$ be a prime element of $I(A,\rho,\BBBB)$
and a prime ideal of $\BBBB$; then $\BBBB\stm I$ is a filter in
$\BBBB$; thus the cluster $\s$ generated by $I$ (see
\ref{primdidclu}) has the property that $\s\cap\BBBB$ is a filter.
Let $X=\Psi^a(A,\rho,\BBBB)$ and suppose that $A$ is complete.
Then $\lag(\BBBB)=CR(X)$ and $\s\in X$. Let $F,G\in CR(X)$ and
$\s\in F\cap G$. There exist $a,b\in\BBBB$ such that $F=\lag(a)$
and $G=\lag(b)$. Thus $a,b\in\s\cap\BBBB$. Therefore $a\wedge
b\in\s$. Then $\s\in\lag(a\we b)=\lag(a)\we\lag(b)=F\we G$. Hence,
 $\int(F\cap G)\nes$. So, if $F,G\in CR(X)$ and $\s\in F\cap G$ then
 $\int(F\cap G)\nes$.
  Obviously, the points of the real
line $\mathbb{R}$ with its natural topology have not this
property. Thus the CLCA $(RC(\mathbb{R}),\rho_{\mathbb{R}},
CR(\mathbb{R}))$ is such that no one prime element of the frame
$I(RC(\mathbb{R}),\rho_{\mathbb{R}}, CR(\mathbb{R}))$ is  a prime
ideal of  $CR(\mathbb{R})$.
\end{rem}

\begin{notas}\label{imnota}
\rm Let $(A,\rho,\BBBB)$ and $(B,\eta,\BBBB\ap)$ be LCAs,
$\p:A\lra B$ be a function and $\s\ap$ be a cluster in
$(B,\eta,\BBBB\ap)$. Then we set:
\begin{itemize}
\item $S_{\s\ap}=\{a\in\BBBB\st (\fa b\in A)[(a\llx
b)\rightarrow(\p(b)\in\s\ap)]\}$;

\item $V_{\s\ap}=\{a\in\BBBB\st (\ex b\in S_{\s\ap})(b\llx a)\}$;

\item $J_{\s\ap}=\BBBB\stm S_{\s\ap}$.
\end{itemize}
\end{notas}

\begin{fact}\label{fss}
Let $(A,\rho,\BBBB)$ and $(B,\eta,\BBBB\ap)$ be LCAs and $\p:A\lra
B$ be a function satisfying conditions (DLC1)-(DLC3). Then, for
every cluster $\s\ap$ in $(B,\eta,\BBBB\ap)$, $\ssi=
\{a\in\BBBB\st (\fa b\in A)[(a\llx
b)\rightarrow((\p(b^*))^*\in\s\ap)]\}$.
\end{fact}

\doc Let $a\in\ssi$, $b\in A$ and $a\llx b$. Then $\p(b)\in\s\ap$.
Since, by Lemma \ref{pf1}(b), $\p(b)\le(\p(b^*))^*$, we get that
$(\p(b^*))^*\in\s\ap$. Thus $a\in R=\{a\in\BBBB\st (\fa b\in
A)[(a\llx b)\rightarrow((\p(b^*))^*\in\s\ap)]\}$.   Conversely,
let $a\in R$, $b\in A$ and $a\llx b$. Then, by (BC1), there exists
$c\in\BBBB$ such that $a\llx c\llx b$. Since $a\in R$ and $a\llx
c$, we get that $(\p(c^*))^*\in\s\ap$. Further, by (DLC3),
$(\p(c^*))^*\lle\p(b)$. Hence, $\p(b)\in\s\ap$. Therefore,
$a\in\ssi$. So, $\ssi=R$. \sqs

\begin{lm}\label{jsdi}
Let $(A,\rho,\BBBB)$ and $(B,\eta,\BBBB\ap)$ be LCAs, $\p:A\lra B$
be a function satisfying conditions (DLC1)-(DLC3) (or conditions
(DLC1), (DLC2), (LC3)), and $\s\ap$ be a cluster in
$(B,\eta,\BBBB\ap)$. Then $\jsi$ is a $\d$-ideal of
$(A,\rho,\BBBB)$. If $\s\ap$ is a bounded cluster in
$(B,\eta,\BBBB\ap)$ and $\p$ satisfies, in addition, condition
(DLC4), then $\jsi$ is a prime element of the frame
$I(A,\rho,\BBBB)$ (see \ref{lideal} for this notation).
\end{lm}

\doc  Obviously,
$\jsi=\{a\in\BBBB\st (\ex b\in A)[(a\llx
b)\we(\p(b)\nin\s\ap)]\}$. Since $0\llx 0$, (DLC1) implies that
$0\in\jsi$. It is clear that $\jsi$ is a lower set. Let
$a_1,a_2\in\jsi$. Then, for $i=1,2$, there exists $b_i\in A$ such
that $a_i\llx b_i$ and $\p(b_i)\nin\s\ap$. Since, for $i=1,2$,
$a_i\in\BBBB$, there exists $c_i\in\BBBB$ such that $a_i\llx
c_i\llx b_i$ (by condition (BC1) in Definition \ref{locono}); then
$\p(c_i)\nin\s\ap$. Set $c=c_1\vee c_2$. Now, by Proposition
\ref{axlc3} (resp., by (LC3)), $\p(c)=\p(c_1\vee
c_2)\lle\p(b_1)\vee\p(b_2)$. Since $\p(b_1)\vee\p(b_2)\nin\s\ap$,
we get that $\p(c)\nin\s\ap$. Therefore, $a_1\vee a_2\in\jsi$. All
this shows that $\jsi$ is an ideal of $A$. Let $a\in\jsi$. Then
there exists $b\in A$ such that $a\llx b$ and $\p(b)\nin\s\ap$.
Using again condition (BC1), we get that there exists $c\in\BBBB$
such that $a\llx c\llx b$. Then $c\in\jsi$ and $a\llx c$. So,
$\jsi$ is a $\d$-ideal of $(A,\rho,\BBBB)$.

Let now $\s\ap\cap\BBBB\ap\nes$ and $\p$ satisfies, in addition,
condition (DLC4). Then $\jsi\neq\BBBB$.
 Indeed, there exists
$b\in\s\ap\cap\BBBB\ap$. Then, by (DLC4), there exists $a\in\BBBB$
such that $b\le\p(a)$;  hence $\p(a)\in\s\ap$. This implies that
$a\in\ssi$. Thus  $\jsi\neq\BBBB$. Let $J_1\cap J_2\sbe \jsi$.
Suppose that there exists $a_i\in J_i\stm\jsi$, $i=1,2$. Since
$J_1, J_2$ are $\d$-ideals, there exists $b_i\in J_i$ such that
$a_i\llx b_i$, $i=1,2$. There exists $c_i\in\BBBB$ such that
$a_i\llx c_i\llx b_i$, $i=1,2$. Since $a_i\in\ssi$, we have that
$\p(c_i)\in\s\ap$, $i=1,2$. By \ref{pf1}(b) (respectively,
\ref{remlc3}), we get that $\p(c_i)\lle\p(b_i)$, $i=1,2$. Now,
Proposition \ref{bbcl1}(a) implies that there exists
$d_i\in\BBBB\ap\cap\s\ap$ such that $d_i\lle\p(b_i)$, $i=1,2$.
Then $d_i\llce\p(b_i)$, $i=1,2$. Thus $(\p(b_i))^*\nin\s\ap$,
$i=1,2$. This implies that $\p(b_1)\we\p(b_2)\in\s\ap$, i.e., by
(DLC2), $\p(b_1\we b_2)\in\s\ap$. We have, however, that $b_1\we
b_2\in \jsi$. Hence, there exists $d\in A$ such that $b_1\we
b_2\llx d$ and $\p(d)\nin\s\ap$. Since $\p(b_1\we b_2)\le\p(d)$
and $\p(b_1\we b_2)\in\s\ap$, we get that $\p(d)\in\s\ap$, a
contradiction. Therefore, $J_1\sbe\jsi$ or $J_2\sbe\jsi$. All this
shows that $\jsi$ is a prime element of the frame
$I(A,\rho,\BBBB)$.
 \sqs

\begin{lm}\label{vsfb}
Let $(A,\rho,\BBBB)$ and $(B,\eta,\BBBB\ap)$ be LCAs, $\p:A\lra B$
be a function satisfying conditions (DLC1)-(DLC4) (or conditions
(DLC1), (DLC2), (LC3), (DLC4)), and $\s\ap$ be a bounded cluster
in $(B,\eta,\BBBB\ap)$. Then $\vsi$ is a filter in $(\BBBB,\le)$.
\end{lm}

\doc  This follows immediately from \ref{jsdi} and \ref{primdfil}.
\sqs

\begin{lm}\label{funsi}
Let $(A,\rho,\BBBB)$ and $(B,\eta,\BBBB\ap)$ be LCAs, $\p:A\lra B$
be a function satisfying conditions (DLC1)-(DLC4) (or conditions
(DLC1), (DLC2), (LC3), (DLC4)), and $\s\ap$ be a bounded cluster
in $(B,\eta,\BBBB\ap)$. Then there exists a unique cluster $\s$ in
$(A,\rho,\BBBB)$ such that $\s\cap\BBBB=\ssi$; moreover,
$\s=\{a\in A\st a\rho\ssi\}$.
\end{lm}

\doc  This follows  from \ref{jsdi} and \ref{primdidclu}.
 \sqs

\begin{nota}\label{nopsi}
Let $(A,\rho,\BBBB)$ and $(B,\eta,\BBBB\ap)$ be LCAs and $\p:A\lra
B$ be a function. We set, for every $a\in A$,
$$D_\p(a)=\bigcup\{I_{\p(b)}\st b\in\BBBB, b\llx a\}$$
(see \ref{dideal} for the notation $I_c$).
\end{nota}

\begin{pro}\label{didle}
Let $(A,\rho,\BBBB)$ and $(B,\eta,\BBBB\ap)$ be LCAs and $\p:A\lra
B$ be a monotone function. Then, for every $a\in A$, $D_\p(a)$ is
a $\d$-ideal of $(B,\eta,\BBBB\ap)$.
\end{pro}

\doc Let $a\in A$. We will prove that
$D_\p(a)=\bigvee\{I_{\p(b)}\st b\in\BBBB, b\llx a\}$, where the
join is taken in the frame $I(B,\eta,\BBBB\ap)$ (see \ref{frlid}).
Then, by \ref{frlid}, $D_\p(a)$ will be a $\d$-ideal.

Set $I=\bigvee\{I_{\p(b)}\st b\in\BBBB, b\llx a\}$. The ideal $I$
is generated by $D_\p(a)$. Hence, $D_\p(a)\sbe I$. Conversely, let
$c\in I$. Then there exists $n\in\mathbb{N}^+$ and, for each
$\ion$, there exist $b_i\in\BBBB$ and $c_i\in\BBBB\ap$ such that
$b_i\llx a$, $c_i\lle\p(b_i)$ and $c=\bigvee\{c_i\st\ion\}$. Set
$b=\bigvee\{b_i\st\ion\}$. Then $b\in\BBBB$, $b\llx a$ and
$c\lle\bigvee\{\p(b_i)\st\ion\}\le\p(b)$. Hence $c\lle\p(b)$, and
since $c\in\BBBB\ap$, we get that $c\in I_{\p(b)}$, where $b\llx
a$.  Thus $c\in D_\p(a)$. \sqs

\begin{lm}\label{funcon}
Let $(A,\rho,\BBBB)$ and $(B,\eta,\BBBB\ap)$ be LCAs, $\p:A\lra B$
be a function satisfying conditions (DLC1)-(DLC4) (or conditions
(DLC1), (DLC2), (LC3), (DLC4)), $X=\Psi^a(A,\rho,\BBBB)$ and
$Y=\Psi^a(B,\eta,\BBBB\ap)$ (see \ref{roeperl} for the notation
$\Psi^a$). For every $\s\ap\in Y$, set $f_\p(\s\ap)=\s$, where
$\s$ is the unique bounded cluster in $(A,\rho,\BBBB)$ such that
$\s\cap\BBBB=\ssi$ (see Lemma \ref{funsi} for $\s$). Then
$f_\p:Y\lra X$ is a continuous function and
\begin{equation}\label{confun}
\fa a\in\BBBB,\  f_\p\inv(\int(\lag(a))=\iota_B(D_\p(a))
\end{equation}
(see \ref{opensetsfr} for $\iota$).
\end{lm}

\doc We will first show that the formula (\ref{confun}) takes
place. So, let $a\in\BBBB$. Since $D_\p(a)=\bigcup\{I_{\p(b)}\st
b\in\BBBB, b\llx a\}$, we get that
%
%\begin{equation}\label{conio}
$$\iota_B(D_\p(a))=\bigcup\{\lbg(c)\st (c\in\BBBB\ap)\we
%\mbox{ and }
(\ex b\in\BBBB)[(b\llx a)\we
%\mbox{ and }
(c\lle\p(b))]\}.$$
%\end{equation}
%

Let $\s\ap\in f_\p\inv(\int(\lag(a))$. Then
$f_\p(\s\ap)=\s\in\int(\lag(a))$. Hence $a^*\nin\s$. Now, by
\ref{cluendcor}, there exists $a_1\in A$ such that $a^*\llcr
a_1^*$ and $a_1^*\nin\s$. We get that $a_1\llcr a$ and
$a_1\in\s\cap\BBBB=\ssi$. Since $a_1\llx a$, there exist
$a_2,b\in\BBBB$ such that $a_1\llx a_2\llx b\llx a$. Then, by the
definition of the set $\ssi$, $\p(a_2)\in\s\ap$. By \ref{pf1}(b)
(resp., by \ref{remlc3}), $\p(a_2)\lle\p(b)$. Now, Proposition
\ref{bbcl1}(a) implies that there exists $c\in\BBBB\ap\cap\s\ap$
such that $c\lle\p(b)$. Thus $\s\ap\in\lbg(c)$, where
$c\in\BBBB\ap$, $c\lle\p(b)$ and $b\llx a$. This means that
$\s\ap\in\iota_B(D_\p(a))$. Hence,
$f_\p\inv(\int(\lag(a))\sbe\iota_B(D_\p(a))$.

Conversely, let $\s\ap\in\iota_B(D_\p(a))$ and $\s=f_\p(\s\ap)$.
Then there exist $b\in\BBBB$ and $c\in\BBBB\ap$ such that $b\llx
a$, $c\lle\p(b)$ and $\s\ap\in\lbg(c)$. Thus $c\in\s\ap$ and hence
$\p(b)\in\s\ap$. This implies that $b\in\ssi=\BBBB\cap\s$. Since
$b\in\BBBB$ and $b\llx a$, we get that $b\llcr a$, i.e.
$b(-C_\rho)a^*$. Thus $a^*\nin\s$. This means that
$f_\p(\s\ap)=\s\in\int(\lag(a))$. Therefore, $\s\ap\in
f_\p\inv(\int(\lag(a))$. We have proved that
$f_\p\inv(\int(\lag(a))\spe\iota_B(D_\p(a))$.

So, the formula (\ref{confun}) is established. Now, by
(\ref{eel}), $\{\int\lag(a)\st a\in\BBBB\}$ is a base of $X$, and,
for every $a\in A$, $D_\p(a)$ is a $\d$-ideal (see Proposition
\ref{didle}). Hence, Theorem \ref{opensetsfr} implies that, for
every $a\in A$, $\iota_B(D_\p(a))$ is an open subset of $Y$. Thus,
by formula (\ref{confun}), $f_\p$ is a continuous function. \sqs

\begin{lm}\label{phief}
Let $f\in \HLC(X,Y)$. Define a function
$\p_f:\Psi^t(Y)\lra\Psi^t(X)$ by the formula:
\begin{equation}\label{xitn}
\fa\, G\in RC(Y), \ \p_f(G)=cl_X(f\inv(\int_Y(G)))
\end{equation}
(see Theorem \ref{roeperl} for $\Psi^t$). Then the function $\p_f$
satisfies conditions (DLC1)-(DLC5) from Definition \ref{dhc} and,
moreover, it satisfies conditions (DLC3S) and (LC3S).
\end{lm}

\doc
Obviously, condition (DLC1) is fulfilled. For proving  condition
(DLC2), recall that (see \cite{dV}) if $U$ and $V$ are two open
subsets of a topological space $Z$ then
\begin{equation}\label{intdv}
\int(\cl(U\cap V))=\int(\cl(U)\cap\cl(V)).
\end{equation}
Let $F,G\in RC(Y)$.  Using  the fact that $\int(F\cap G)$ is a
regular open set, we get that $\int(F\cap G)=\int(\cl(\int(F\cap
G)))$. Thus
$$\p_f(F\we G)=\cl(f\inv(\int(\cl(\int(F\cap
G)))))=\cl(f\inv(\int(F\cap G))).$$
Now, setting $U=f\inv(\int(F))$ and $V=f\inv(\int(G))$, we obtain,
using (\ref{intdv}), that
$$\p_f(F)\we\p_f(G)=\cl(U)\we\cl(V)=
\cl(\int(\cl(U)\cap\cl(V)))=$$
$$=\cl(\int(\cl(U\cap V)))=\cl(U\cap V)=\cl(f\inv(\int(F\cap
G))).$$
 Therefore, $\p_f(F\we G)=\p_f(F)\we\p_f(G)$. So, (DLC2)
is fulfilled.

We will now show that not only condition (DLC3) is true, but even
condition (DLC3S) takes place.  Indeed, let $F,G\in RC(Y)$ and
$F\ll_{\rho_Y} G$. Then $F\sbe\int(G)$ and
$(\p_f(F^*))^*=(\cl(f\inv(\int(F^*))))^*=(\cl(f\inv(Y\stm
F)))^*=(\cl(X\stm f\inv(F)))^*=\cl(\int(f\inv(F)))\sbe
f\inv(F)\sbe
f\inv(\int(G))\sbe\int(\cl(f\inv(\int(G))))=\int(\p_f(G))$. Hence
$(\p_f(F^*))^*\ll_{\rho_X} \p_f(G)$, i.e. condition (DLC3S) is
fulfilled.

For verifying (DLC4), let $H\in CR(X)$. Then $f(H)$ is compact.
Since $Y$ is locally compact, there exists $F\in CR(Y)$ such that
$f(H)\sbe\int(F)$. Now we obtain that $H\sbe
f\inv(\int(F))\sbe\int(\cl(f\inv(\int(F))))=\int(\p_f(F))$, i.e.
$H\ll_{\rho_X}\p_f(F)$. Hence, condition (DLC4) takes place.

Let $F\in RC(Y)$. For establishing condition (DLC5), we have to
show that $\p_f(F)=\bigvee\{\p_f(G)\st G\in CR(Y),
G\sbe\int(F)\}$. We have that $\p_f(F)=\cl(f\inv(\int(F)))$ and,
since $Y$ is locally compact and regular,
$\int(F)=\bigcup\{\int(G)\st G\in CR(Y), G\sbe\int(F)\}$. Further,
it is obvious that if
 $G\sbe\int(F)$ then $\cl(f\inv(\int(G)))\sbe f\inv(G)\sbe
f\inv(\int(F))$. Now, it is easy to see that the desired equality
is fulfilled. So, condition (DLC5) is verified.

Finally, we will show that condition (LC3S) is fulfilled as well.
Let, for $i=1,2$, $F_i,G_i\in RC(Y)$ and $F_i\sbe\int(G_i)$. We
have to show that $\p_f(F_1\cup
F_2)\sbe\int(\p_f(G_1)\cup\p_f(G_2))$. Indeed, $\p_f(F_1\cup
F_2)=\cl(f\inv(\int(F_1\cup F_2)))\sbe f\inv(F_1\cup F_2)\sbe
f\inv(\int(G_1)\cup\int(G_2))\sbe\int(\cl(f\inv(\int(G_1)\cup\int(G_2))))=
\int(\p_f(G_1)\cup\p_f(G_2))$. So, condition (LC3S) is verified.
\sqs

\begin{lm}\label{commdiagr}
Let $(A,\rho,\BBBB)$ and $(B,\eta,\BBBB\ap)$ be CLCAs, $\p:A\lra
B$ be a function satisfying conditions (DLC1)-(DLC5) (or
conditions (DLC1), (DLC2), (LC3), (DLC4), (DLC5)),
$X=\Psi^a(A,\rho,\BBBB)$ and $Y=\Psi^a(B,\eta,\BBBB\ap)$ (see
\ref{roeperl} for the notation $\Psi^a$). Let $f=f_\p$ (see Lemma
\ref{funcon} for $f_\p$) and $\p\ap=\p_f$ (see Lemma \ref{phief}
for $\p_f$). Then $\lbg\circ\p=\p\ap\circ\lag$ (see \ref{roeperl}
for the notations $\lag$ and $\lbg$).
\end{lm}

\doc Note that, by Lemma \ref{funcon}, $f:Y\lra X$ is a continuous
function. Hence, Lemma \ref{phief} implies that the function
$\p\ap$ satisfies conditions (DLC1)-(DLC5).

Let us now regard the case when $a\in \BBBB$. We have to show that
$\lbg(\p(a))=\p\ap(\lag(a))$. By the definitions of $\p\ap$ and
$f$, and the formula (\ref{confun}), we obtain that
$\p\ap(\lag(a))=\cl(f_\p\inv(\int(\lag(a))))=\cl(\iota_B(D_\p(a)))=
\cl(\bigcup\{\lbg(b)\st b\in D_\p(a)\}=\bigvee\{\lbg(b)\st b\in
D_\p(a)\}=\lbg(\bigvee\{b\st b\in D_\p(a)\})$ (since, by Theorem
\ref{roeperl}, $\lbg$ is an LCA-isomorphism). Hence, we have to
prove that $\p(a)=\bigvee\{b\st b\in D_\p(a)\})=\bigvee D_\p(a)$.
By condition (DLC5), we have that $\p(a)=\bigvee\{\p(c)\st
c\in\BBBB, c\llx a\}$. Since, for every $c\in A$,
$\p(c)=\bigvee\{b\in\BBBB\ap\st b\lle\p(c)\}$, we get that
$\p(a)=\bigvee\{b\in\BBBB\ap\st (\ex c\in\BBBB)[(c\llx a)\we
(b\lle\p(c))]\}$. By definition, $D_\p(a)=\bigcup\{I_{\p(c)}\st
c\in\BBBB, c\llx a\}$. Thus $(b\in
D_\p(a))\leftrightarrow[(b\in\BBBB\ap)\we ((\ex c\in\BBBB)((c\llx
a)\we (b\lle\p(c))))]$. This shows that $\p(a)=\bigvee D_\p(a)$.
So, we have proved that $\lbg(\p(a))=\p\ap(\lag(a))$ for every
$a\in\BBBB$.

Let now $a\in A$. Then, by condition (DLC5),
$\p(a)=\bigvee\{\p(b)\st b\in\BBBB, b\llx a\}$. Hence, using the
fact that $\lbg$ is an LCA-isomorphism and the formula proved in
the preceding paragraph, we get that
$\lbg(\p(a))=\bigvee\{\lbg(\p(b))\st b\in\BBBB, b\llx
a\}=\bigvee\{\p\ap(\lag(b))\st b\in\BBBB, b\llx a\}$. Further,
since the function $\p\ap$ satisfies condition (DLC5), we have
that for every $G\in RC(X)$, $\p\ap(G)=\bigvee\{\p\ap(F)\st F\in
CR(X), F\ll_{\rho_X} G\}$. Now using the fact that $\lag$ is an
LCA-isomorphism between LCAs $(A,\rho,\BBBB)$ and
$(RC(X),\rho_X,CR(X))$, we get that
$\p\ap(\lag(a))=\bigvee\{\p\ap(\lag(b))\st b\in\BBBB, b\llx
a\}=\lbg(\p(a))$. So, the desired equality is established. \sqs

\begin{lm}\label{dlc3slc3s}
Let $(A,\rho,\BBBB)$ and $(B,\eta,\BBBB\ap)$ be CLCAs and
$\p:A\lra B$ be a function satisfying conditions (DLC1)-(DLC5) (or
conditions (DLC1), (DLC2), (LC3), (DLC4), (DLC5)). Then $\p$
satisfies conditions (DLC3S) and (LC3S) as well.
\end{lm}

\doc  Let $f=f_\p$ (see Lemma
\ref{funcon} for $f_\p$) and $\p\ap=\p_f$ (see Lemma \ref{phief}
for $\p_f$). Then, by Lemma \ref{commdiagr},
$\lbg\circ\p=\p\ap\circ\lag$. Since the function $\p\ap$ satisfies
conditions (DLC3S) and (LC3S) (by Lemma \ref{phief}) and the
functions $\lag$ and $\lbg$ are LCA-isomorphisms, we get that the
function $\p$ satisfies conditions (DLC3S) and (LC3S) as well.
\sqs

The above lemma implies the following fact mentioned in the
previous section:

\begin{cor}\label{dlc3sax}
Condition (DLC3) in Definition \ref{dhc} can be replaced by any of
the conditions (DLC3S), (LC3) and (LC3S) (i.e., we obtain
equivalent systems of axioms by these replacements).
\end{cor}

\begin{lm}\label{compos1}
Let $\p:(A,\rho,\BBBB)\lra (B,\eta,\BBBB\ap)$ be a function
between CLCAs and let $\p$  satisfy conditions (DLC1)-(DLC4). Then
the function $\p\cuk$ (see (\ref{cukf})) satisfies conditions
(DLC1)-(DLC5).
\end{lm}

\doc Obviously, for every $a\in A$, $\p\cuk(a)\le\p(a)$. Hence,
$\p\cuk(0)=0$, i.e. (DLC1) is fulfilled. For (DLC2) and (DLC5) see
\ref{pf1}(d).

Let $a\in\BBBB, b\in A$ and $a\llx b$. Then, by (BC1), there exist
$c,d\in\BBBB$ such that $a\llx c\llx d\llx b$. Thus $a\llcr c$ and
hence $c^*\llcr a^*$. Now, using Lemma \ref{cuklm},  we obtain
that $\p(c^*)\le\p\cuk(a^*)$. Since $\p(d)\le\p\cuk(b)$, we get
that
  $(\p\cuk(a^*))^*\le (\p(c^*))^*\lle
\p(d)\le\p\cuk(b)$. Therefore, $(\p\cuk(a^*))^*\lle\p\cuk(b)$. So,
(DLC3) is fulfilled.

For verifying (DLC4), let $b\in\BBBB$. Then there exists
$a\in\BBBB$ such that $b\le\p(a)$. By (BC1), there exists
$a_1\in\BBBB$ with $a\llx a_1$. Then $b\le\p(a)\le\p\cuk(a_1)$.
Thus, $\p\cuk$ satisfies condition (DLC4). \sqs

\begin{lm}\label{assocuk}
Let $\p_i:(A_i,\rho_i,\BBBB_i)\lra
(A_{i+1},\rho_{i+1},\BBBB_{i+1})$, where $i=1,2$, be two functions
between CLCAs. Then:

\smallskip

\noindent(a) $(\p_2\cuk\circ\p_1)\cuk=(\p_2\circ\p_1)\cuk$;

\noindent(b) If $\p_1$ and $\p_2$ are monotone functions, then
$(\p_2\circ\p_1\cuk)\cuk=(\p_2\circ\p_1)\cuk$;

\noindent(c) If $\p_1$ and $\p_2$ satisfy conditions (DLC1)-(DLC5)
then the function $\p_2\circ\p_1$ satisfies conditions
(DLC1)-(DLC4) and even condition (DLC3S).
\end{lm}

\doc   We will write, for $i=1,2$,
$``\ll_i$" instead of $``\ll_{\rho_i}$". We also set
$\p=\p_2\circ\p_1$.

\smallskip

\noindent(a)  Let $a\in A_1$. Then
$(\p_2\cuk\circ\p_1)\cuk(a)=\bv\{\p_2\cuk(\p_1(b))\st b\in\BBBB_1,
b\ll_1 a\}=\bv\{\bv\{\p(c)\st c\in\BBBB_1, c\ll_1 b\}\st
b\in\BBBB_1, b\ll_1 a\}=\bv\{\p(c)\st c\in\BBBB_1, c\ll_1
a\}=\p\cuk(a)$.

\smallskip

\noindent(b) Let $a\in A_1$. Then $L=(\p_2\circ\p_1\cuk)\cuk(a)=
\bv\{\p_2(\p_1\cuk(b))\st b\in\BBBB_1, b\ll_1
a\}=\bv\{\p_2(\bv\{\p_1(c)\st c\in\BBBB_1, c\ll_1 b\})\st
b\in\BBBB_1, b\ll_1 a\}$ and $R=\p\cuk(a)=\bv\{\p(c)\st
c\in\BBBB_1, c\ll_1 a\}$. Let $c\in\BBBB_1$ and $c\ll_1 a$. Then,
by (BC1), there exists $b\in\BBBB_1$ such that $c\ll_1 b\ll_1 a$.
This shows that $R\le L$. Conversely, let $b\in\BBBB_1$ and
$b\ll_1 a$. Then, for every $c\in\BBBB_1$ such that $c\ll_1 b$, we
have that $\p_1(c)\le\p_1(b)$. Hence $\bv\{\p_1(c)\st c\in\BBBB_1,
c\ll_1 b\}\le\p_1(b)$. Then $\p_2(\bv\{\p_1(c)\st c\in\BBBB_1,
c\ll_1 b\})\le\p(b)$ and $b\ll_1 a$. This implies that $L\le R$.
So, $L=R$. Hence, $(\p_2\circ\p_1\cuk)\cuk=(\p_2\circ\p_1)\cuk$.

\smallskip

\noindent(c) Obviously, the function $\p$ satisfies conditions
(DLC1), (DLC2) and (DLC4). For proving that $\p$ satisfies
condition (DLC3S),
 let   $a,b\in A_1$ and $a\ll_1 b$.
Since the functions  $\p_1$ and $\p_2$ satisfy condition
 (DLC3S) (by Lemma \ref{dlc3slc3s}), we obtain that $(\p_1(a^*))^*\ll_2
 \p_1(b)$ and
$(\p_2(\p_1(a^*)))^*\ll_3\p_2(\p_1(b))$, i.e.
$(\p(a^*))^*\ll_3\p(b)$. Hence, the function $\p$ satisfies
condition (DLC3S).
 \sqs

\begin{pro}\label{dlccat}
$\DLC$ is a category.
\end{pro}

\doc It is clear that for every CLCA $(A,\rho,\BBBB)$, the usual
identity function $id_A:A\lra A$ satisfies conditions
(DLC1)-(DLC5); moreover, using Lemma \ref{pf1}(e), we get that if
$(B,\eta,\BBBB\ap)$ and $(B_1,\eta_1,\BBBB\ap_1)$ are CLCAs, and
$\p:(A,\rho,\BBBB)\lra(B,\eta,\BBBB\ap)$ and
$\psi:(B_1,\eta_1,\BBBB\ap_1)\lra(A,\rho,\BBBB)$ are functions
satisfying condition (DLC5), then $id_A\cd\psi=\psi$ and $\p\cd
id_A=\p$. So, $id_A$ is the $\DLC$-identity on $(A,\rho,\BBBB)$.

Let $\p_i:(A_i,\rho_i,\BBBB_i)\lra
(A_{i+1},\rho_{i+1},\BBBB_{i+1})$, where $i=1,2$, be two functions
between CLCAs, and let $\p_1$ and $\p_2$ satisfy conditions
(DLC1)-(DLC5). We will show that  the function $\p_2\cd\p_1$
satisfies conditions (DLC1)-(DLC5).

Set $\p=\p_2\circ\p_1$.
 Then, by Lemma \ref{assocuk}(c),  the function $\p$
satisfies conditions (DLC1)-(DLC4). Now,  Lemma \ref{compos1}
implies that the function $\p\cuk$ satisfies conditions
(DLC1)-(DLC5). Since $\p_2\cd\p_1=\p\cuk$, we get that the
function $\p_2\cd\p_1$ satisfies conditions (DLC1)-(DLC5).

Finally, we will show that the composition in $\DLC$ is
associative. Let, for $i=1,2,3$, $\p_i:(A_i,\rho_i,\BBBB_i)\lra
(A_{i+1},\rho_{i+1},\BBBB_{i+1})$  be a function between CLCAs
satisfying conditions (DLC1)-(DLC5). We will show that
$(\p_3\cd\p_2)\cd\p_1=\p_3\cd(\p_2\cd\p_1)$. Using Lemma
\ref{assocuk}, we get that
$(\p_3\cd\p_2)\cd\p_1=((\p_3\circ\p_2)\cuk\circ\p_1)\cuk=((\p_3\circ\p_2)\circ\p_1)\cuk$
and
$\p_3\cd(\p_2\cd\p_1)=(\p_3\circ(\p_2\circ\p_1)\cuk)\cuk=(\p_3\circ(\p_2\circ\p_1))\cuk$.
Thus, the associativity of the composition in $\DLC$ is proved.

All this shows that $\DLC$ is a category. \sqs

\begin{lm}\label{fcuk}
Let $\p:(A,\rho,\BBBB)\lra(B,\eta,\BBBB\ap)$ be a function between
CLCAs satisfying conditions (DLC1)-(DLC4). Then $f_\p=f_{\p\cuk}$
(see Lemma \ref{funcon} for the notation $f_\p$ and
(\ref{cukfcon}) for the notation $\p\cuk$).
\end{lm}

\doc Let $X=\Psi^a(A,\rho,\BBBB)$ and $Y=\Psi^a(B,\eta,\BBBB\ap)$.
By Lemma \ref{compos1}, the function $\p\cuk$ satisfies conditions
(DLC1)-(DLC4). Hence, we can apply
 Lemma \ref{funcon} in order to construct two (continuous) functions
  $f_\p,f_{\p\cuk}:Y\lra X$.  Let $\s\ap\in
Y$. Set $\s=f_\p(\s\ap)$ and $\s_1=f_{\p\cuk}(\s\ap)$. By
\ref{bbcl}, for proving that $\s=\s_1$, it is enough to show that
$\s\cap\BBBB=\s_1\cap\BBBB$, where $\s\cap\BBBB=\{a\in\BBBB\st
(\fa b\in A)[ (a\llx b)\rightarrow (\p(b)\in\s\ap)]\}$
 and $\s_1\cap\BBBB=\{a\in\BBBB\st (\fa
b\in A)[ (a\llx b)\rightarrow (\p\cuk(b)\in\s\ap)]\}$   (see Lemma
\ref{funcon}).

Let $a\in\s_1\cap\BBBB$, $b\in A$ and $a\llx b$. Then
$\p\cuk(b)\in\s\ap$. Since $\p\cuk(b)\le\p(b)$ (by \ref{pf1}(g)),
we get that $\p(b)\in\s\ap$. So, $\s\cap\BBBB\spe\s_1\cap\BBBB$.
Conversely, let $a\in\s\cap\BBBB$, $b\in A$ and $a\llx b$. By
(BC1), there exists $c\in\BBBB$ such that $a\llx c\llx b$. Then
$\p(c)\in\s\ap$ and $\p(c)\le\p\cuk(b)$. Hence,
$\p\cuk(b)\in\s\ap$. So, $\s\cap\BBBB\sbe\s_1\cap\BBBB$.
Therefore, $\s=\s_1$. This shows that $f_\p=f_{\p\cuk}$. \sqs

\begin{pro}\label{funla}
For every $(A,\rho,\BBBB)\in\card\DLC$, set
$$\LAM^a(A,\rho,\BBBB)=\Psi^a(A,\rho,\BBBB)$$
(see \ref{roeperl} for $\Psi^a$), and for every
$\p\in\DLC((A,\rho,\BBBB),(B,\eta,\BBBB\ap))$, define
$$\LAM^a(\p):\LAM^a(B,\eta,\BBBB\ap)\lra\LAM^a(A,\rho,\BBBB)$$
by the formula $\LAM^a(\p)=f_\p$, where $f_\p$ is the function
defined in Lemma \ref{funcon}.
 Then $\Lambda^a:\DLC\lra\HLC$ is a contravariant functor.
\end{pro}

\doc By Theorem \ref{roeperl}, if
$(A,\rho,\BBBB)\in\card\DLC$ then
$\LAM^a(A,\rho,\BBBB)\in\card\HLC$, and, by  Lemma \ref{funcon},
 if
$\p\in\DLC((A,\rho,\BBBB),(B,\eta,\BBBB\ap))$ then
$\LAM^a(\p)\in\HLC(\LAM^a(B,\eta,\BBBB\ap),\LAM^a(A,\rho,\BBBB))$.
Further, let $(A,\rho,\BBBB)\in\card\DLC$ and set
$X=\LAM^a(A,\rho,\BBBB)$, $f=\LAM^a(id_A)$. We have to show that
$f=id_X$. Indeed, let $\s\ap\in X$. Set $\s=f(\s\ap)$. We will
prove that $\s\ap\cap\BBBB=\s\cap\BBBB$; then Corollary \ref{bbcl}
will imply that $\s=\s\ap$. We have, by the definition of $\s$,
that $\s\cap\BBBB=\{a\in\BBBB\st (\fa b\in A)[ (a\llx
b)\rightarrow (b\in\s\ap)]\}$. Obviously,
$\s\ap\cap\BBBB\sbe\s\cap\BBBB$. Conversely, let
$a\in\BBBB\cap\s$. Suppose that $a\nin\s\ap$. Then, by
\ref{cluendcor}, there exists $b\in A$ such that $a\llcr b$ and
$b\nin\s\ap$. Since $a\llx b$, we have that $b\in\s\ap$, a
contradiction. Hence, $\s\ap\cap\BBBB\spe\s\cap\BBBB$. So,
$f=id_X$.

Let $\p_i\in\DLC((A_i,\rho_i,\BBBB_i),
(A_{i+1},\rho_{i+1},\BBBB_{i+1}))$, where $i=1,2$. Set, for
$i=1,2,3$, $X_i=\LAM^a(A_i,\rho_i,\BBBB_i)$, and,  for $i=1,2$,
$f_i=\LAM^a(\p_i)$. We will write, for $i=1,2,3$, $``\ll_i$"
instead of $``\ll_{\rho_i}$". Let $\p=\p_2\circ\p_1$ and
$f=f_1\circ f_2$.
 We have to show that $\LAM^a(\p_2\cd\p_1)=f$.
By Lemma \ref{assocuk}(c), the function $\p$ satisfies conditions
(DLC1)-(DLC4). Thus, by Lemma \ref{funcon},  the function
$g=f_\p:X_3\lra X_1$ is well-defined. We will show that $g=f$. Let
$\s_3\in X_3$ and $\s=g(\s_3)$. Then
$\s\cap\BBBB_1=\{a\in\BBBB_1\st (\fa b\in A_1) [(a\ll_1
b)\rightarrow(\p(b)\in\s_3)]\}$. Let $\s_2=f_2(\s_3)$ and
$\s_1=f_1(\s_2)$.  For proving that $\s=\s_1$, it is enough to
show (by \ref{bbcl}) that $\s\cap\BBBB_1=\s_1\cap\BBBB_1$. We have
that $\s_2\cap\BBBB_2=\{a\in\BBBB_2\st (\fa b\in A_2) [ (a\ll_2
b)\rightarrow(\p_2(b)\in\s_3)]\}= \{a\in\BBBB_2\st (\fa b\in A_2)
[ (a\ll_2 b)\rightarrow((\p_2(b^*))^*\in\s_3)\}$ and
 $\s_1\cap\BBBB_1=\{a\in\BBBB_1\st (\fa b\in A_1) [ (a\ll_1
b)\rightarrow(\p_1(b)\in\s_2)]\}=\{a\in\BBBB_1\st (\fa b\in A_1) [
(a\ll_1 b)\rightarrow((\p_1(b^*))^*\in\s_2)]\}$ (see \ref{funcon}
and \ref{fss}). Let $a\in\s_1\cap\BBBB_1$, $b\in A_1$ and $a\ll_1
b$. Then, by (BC1), there exists $c\in\BBBB_1$ such that $a\ll_1
c\ll_1 b$. Then $\p_1(c)\in\s_2$ and, by \ref{pf1}(b),
$\p_1(c)\ll_2\p_1(b)$. Hence, by Proposition \ref{bbcl1}(a), there
exists $d_2\in\BBBB_2\cap\s_2$ such that $d_2\ll_2\p_1(b)$. Then
$\p(b)\in\s_3$. So, $a\in\s\cap\BBBB_1$. Thus,
$\s_1\cap\BBBB_1\sbe\s\cap\BBBB_1$. Conversely, let
$a\in\s\cap\BBBB_1$. Suppose that $a\nin\s_1\cap\BBBB_1$. Then
there exists $b\in A_1$ such that $a\ll_1 b$ and
$(\p_1(b^*))^*\nin\s_2$. There exists $c\in\BBBB_1$ such that
$a\ll_1 c\ll_1 b$. Since $\p_1(b^*)\in\s_2$ and, by (DLC3S),
$\p_1(b^*)\ll_2\p_1(c^*)$, Proposition \ref{bbcl1}(a) implies that
there exists $d_1\in\s_2\cap\BBBB_2$ such that
$d_1\ll_2\p_1(c^*)$. Thus $\p(c^*)\in\s_3$. Further, there exists
$d\in\BBBB_1$ such that $a\ll_1 d\ll_1 c$. Then $\p(d)\in\s_3$
and, by (DLC3S) (see \ref{assocuk}(c)), $\p(d)\ll_3 \p(c)$. Using
again Proposition \ref{bbcl1}(a), we get that there exists
$e\in\s_3\cap\BBBB_3$ such that $e\ll_3\p(c)$. Thus
$e\ll_{C_{\rho_3}}\p(c)$, i.e. $e(-C_{\rho_3})(\p(c))^*$. Then, by
\ref{pf1}(b), $e(-C_{\rho_3})\p(c^*)$. Therefore
$\p(c^*)\nin\s_3$, a contradiction. It shows that
$a\in\s_1\cap\BBBB_1$. So, $\s_1\cap\BBBB_1\spe\s\cap\BBBB_1$. We
have proved that $\s=\s_1$. So, $g=f$. Since $\p_2\cd\p_1=\p\cuk$,
Lemma \ref{fcuk} implies that  $\LAM^a(\p_2\cd\p_1)=g$. Therefore,
 $\LAM^a(\p_2\cd\p_1)=f$. \sqs

\begin{pro}\label{funlt}
For every $X\in\card\HLC$, set
$\LAM^t(X)=\Psi^t(X)$
(see \ref{roeperl} for $\Psi^t$), and for every $f\in\HLC(X,Y)$,
define
$\LAM^t(f):\LAM^t(Y)\lra\LAM^t(X)$
by the formula $\LAM^t(f)=\p_f$, where $\p_f$ is the function
defined in Lemma \ref{phief}.
 Then $\Lambda^t:\HLC\lra\DLC$ is a contravariant functor.
\end{pro}

\doc  By Theorem \ref{roeperl}, if
$X\in\card\HLC$ then $\LAM^t(X)\in\card\DLC$, and, by Lemma
\ref{phief},
 if
$f\in\HLC(X,Y)$ then $\LAM^t(f)\in\DLC(\LAM^t(Y),\LAM^t(X))$.
Further, it is obvious that $\LAM^t$ preserves identity morphisms.

Let $f\in\HLC(X,Y)$ and $g\in\HLC(Y,Z)$. We will prove that
$\LAM^t(g\circ f)=\LAM^t(f)\diamond\LAM^t(g)$. Set $h=g\circ f$.
We have that $\p_h=\LAM^t(h)$, $\p_f=\LAM^t(f)$ and
$\p_g=\LAM^t(g)$ (see  Lemma \ref{phief} for the notations $\p_f$
etc.). Let $F\in RC(Z)$. Then
$\p_h(F)=\cl(h\inv(\int(F)))=\cl(f\inv(g\inv(\int(F))))$ and
$(\p_f\circ\p_g)\cuk(F)=\bigvee\{\p_f(\p_g(G))\st G\in CR(Z),\
G\ll_{\rho_Z} F\}$. If $G\in CR(Z)$ (or even $G\in RC(Z)$) and
$G\sbe\int(F)$, then
$$\p_f(\p_g(G))\sbe f\inv(\cl(g\inv(\int(G))))\sbe
f\inv(g\inv(G))\sbe f\inv(g\inv(\int(F))).$$
 Thus
$(\p_f\circ\p_g)\cuk(F)\sbe\p_h(F)$. Further, since
$\int(F)=\bigcup\{\int(G)\st G\in CR(Z), G\sbe\int(F)\}$, we get
that $g\inv(\int(F))\sbe \bigcup\{\int(\p_g(G))\st G\in CR(Z),
G\sbe\int(F)\}$. Hence
$$\p_h(F)\sbe\cl(\bigcup\{f\inv(\int(\p_g(G)))\st G\in CR(Z),
G\sbe\int(F)\})\sbe$$
$$\sbe\cl(\bigcup\{\p_f(\p_g(G))\st G\in CR(Z),
G\sbe\int(F)\})=(\p_f\circ\p_g)\cuk(F).$$
Therefore, $\p_h=\p_f\diamond\p_g$. So, $\Lambda^t:\HLC\lra\DLC$
is a contravariant functor. \sqs

\begin{pro}\label{natisomta}
The identity functor $Id_{\,\DLC}$ and the  functor
$\LAM^t\circ\LAM^a$
 are naturally isomorphic.
\end{pro}

\doc Recall that for every $(A,\rho,\BBBB)\in\card{\DLC}$, the function
$$\l_A^g:(A,\rho,\BBBB)\lra(\LAM^t\circ\LAM^a)(A,\rho,\BBBB)$$
 is an
LCA-isomorphism (see (\ref{hapisomn})). We will show that
$\l^g:Id_{\,\DLC}\lra \LAM^t\circ\LAM^a$, where for every
$(A,\rho,\BBBB)\in\card\DLC$, $\l^g(A,\rho,\BBBB)=\l_A^g$, is a
natural isomorphism. (Note that, clearly, every LCA-isomorphism is
a $\DLC$-isomorphism.)

Let $\p\in\DLC((A,\rho,\BBBB),(B,\eta,\BBBB\ap))$. We have to show
that $\lbg\cd\p=(\LAM^t\circ\LAM^a)(\p)\cd\lag$, i.e. that
$(\lbg\circ\p)\cuk=((\LAM^t\circ\LAM^a)(\p)\circ\lag)\cuk$. Since,
by Lemma \ref{commdiagr} and the definitions of the contravariant
functors $\LAM^t$ and $\LAM^a$, we have that
\begin{equation}\label{labggg}
\lbg\circ\p=(\LAM^t\circ\LAM^a)(\p)\circ\lag,
\end{equation}
our assertion follows immediately. So, $\l^g$ is a natural
isomorphism. \sqs

\begin{pro}\label{natisomat}
The identity functor $Id_{\,\HLC}$ and the  functor
$\LAM^a\circ\LAM^t$
 are naturally isomorphic.
\end{pro}

\doc Recall that, for every $X\in\card{\HLC}$, the map
$t_X:X\lra(\LAM^a\circ\LAM^t)(X)$, where $t_X(x)=\s_x$ for every
$x\in X$, is a homeomorphism (see (\ref{homeo}) for $t_X$ and
(\ref{sxvx}) for the notation $\s_x$). We will show that
$t^l:Id_{\,\HLC}\lra \LAM^a\circ\LAM^t$, where for every
$X\in\card\HLC$, $t^l(X)=t_X$, is a natural isomorphism.

Let $f\in\HLC(X,Y)$ and $f\ap=(\LAM^a\circ\LAM^t)(f)$.
%$X\ap=(\LAM^a\circ\LAM^t)(X)$, $Y\ap=(\LAM^a\circ\LAM^t)(Y)$.
We have to prove that $t_Y\circ f=f\ap\circ t_X$, i.e. that for
every $x\in X$, $\s_{f(x)}=f\ap(\s_x)$. By Corollary \ref{bbcl},
it is enough to show that $\s_{f(x)}\cap CR(Y)=f\ap(\s_x)\cap
CR(Y)$.

 We have, by the definition of $\LAM^t$, that
$\LAM^t(f)=\p_f$, where, for every $G\in RC(Y)$,
$\p_f(G)=\cl(f\inv(\int(G)))$. Hence $f\ap(\s_x)\cap CR(Y)=\{F\in
CR(Y)\st (\fa G\in
RC(Y))[(F\sbe\int(G))\rightarrow(\p_f(G)\in\s_x)]\}=\{F\in
CR(Y)\st (\fa G\in
RC(Y))[(F\sbe\int(G))\rightarrow(x\in\cl(f\inv(\int(G))))]\}$. Let
$F\in\s_{f(x)}\cap CR(Y)$. Then $f(x)\in F$. If $G\in RC(Y)$ and
$F\sbe\int(G)$, then $f(x)\in\int(G)$. Thus $x\in
f\inv(\int(G))\sbe\cl(f\inv(\int(G)))$. Therefore $F\in
f\ap(\s_x)\cap CR(Y)$. So,  $\s_{f(x)}\cap CR(Y)\sbe
f\ap(\s_x)\cap CR(Y)$. Conversely, let $F\in f\ap(\s_x)\cap
CR(Y)$. Suppose that $f(x)\nin F$. Then, by \ref{loccome}, there
exists $G\in CR(Y)$ such that $F\sbe\int(G)\sbe G\sbe
Y\stm\{f(x)\}$. Thus $x\nin f\inv(G)$. Since
$\cl(f\inv(\int(G)))\sbe f\inv(G)$, we get that
$x\nin\cl(f\inv(\int(G)))$, a contradiction. Hence $f(x)\in F$,
i.e. $F\in\s_{f(x)}$. Thus $\s_{f(x)}\cap CR(Y)\spe f\ap(\s_x)\cap
CR(Y)$. Therefore, for every $x\in X$, $\s_{f(x)}=f\ap(\s_x)$.
This means that $t_Y\circ f=f\ap\circ t_X$. So, $t^l$ is a natural
isomorphism. \sqs

It is clear now that  Theorem \ref{lccont} follows from
Propositions \ref{dlccat}, \ref{funla}, \ref{funlt},
\ref{natisomta} and \ref{natisomat}. So, the proof of Theorem
\ref{lccont} is complete.

%--------------------------------------------------------------------

\baselineskip = 1.2\normalbaselineskip

\baselineskip = 1.00\normalbaselineskip

\vspace{0.25cm}

 Faculty  of Mathematics and Informatics,

Sofia University,

5 J. Bourchier Blvd.,

1164 Sofia,

Bulgaria

\end{document}